\documentclass[
]{ceurart}

\usepackage{amsfonts}

\usepackage{latexsym}
\usepackage{amssymb}
\usepackage{amsmath}
\usepackage{amsthm}
\usepackage[mathscr]{eucal}
\usepackage{graphicx}
\usepackage{hyperref}
\usepackage{caption}
\usepackage{subcaption}

\newtheorem{Theorem}{Theorem}
\newtheorem{Definition}{Definition}
\newtheorem{Proposition}{Proposition}

\begin{document}

\copyrightyear{2021}
\copyrightclause{Copyright for this paper by its authors.
  Use permitted under Creative Commons License Attribution 4.0
  International (CC BY 4.0).}

\conference{Algorithms, Computing and Mathematics Conference (ACM 2021), Chennai, India}

\title{Duality and Outermost Boundaries in Generalized Percolation Lattices}

\author[1]{Ghurumuruhan Ganesan}[%
email=gganesan82@gmail.com,
]
\address[1]{Institute of Mathematical Sciences, HBNI, Chennai}

\begin{abstract}
In this paper we consider a connected planar graph~\(G\) and impose conditions
that results in~\(G\) having a percolation lattice-like cellular structure. Assigning
each cell of~\(G\) to be either occupied or vacant, we describe the outermost boundaries of star and plus connected components in~\(G.\) We then consider the dual graph of~\(G\) and impose conditions under which the dual is also a percolation lattice. Finally, using~\(G\) and its dual, we construct vacant cell cycles surrounding occupied components and study left right crossings and bond percolation in rectangles.
\end{abstract}


\maketitle

\renewcommand{\theequation}{\thesection.\arabic{equation}}
\setcounter{equation}{0}
\section{Introduction} \label{intro}
The structure of the outermost boundary of finite components is crucial for contour analysis problems of percolation~\cite{grim} and random graphs~\cite{penrose}. For the square lattice, self-duality plays a crucial role in determining the properties of star and plus connected components and we refer to Chapter~\(3\)~\cite{boll2} for a detailed discussion of combinatorial properties of percolation in regular lattices. For general graphs,~\cite{timar} uses separating sets in equivalence class of infinite paths to study duality in locally finite graphs and~\cite{penrose} uses unicoherence and topological arguments to investigate plus connected components.

In many applications, it might happen that the lattice on which percolation occurs is not necessarily regular, like for example percolation in Voronoi tessellations~\cite{boll2}. It would therefore be interesting to study the duality properties of such irregular lattices and determine conditions under these lattices have behaviour  similar to the regular lattices. In this paper we study outermost boundaries in generalized percolation lattices and prove duality properties analogous to regular lattices. We first consider an arbitrary planar graph~\(G\) and impose certain cyclic conditions on~\(G\) that results in a cellular structure analogous to regular lattices. We then define the dual graph of~\(G\) and determine necessary and sufficient conditions for the dual to be a percolation lattice, analogous to~\(G.\) Using~\(G\) and its dual, we study outermost boundaries, occupied components and rectangular left-right crossings in generalized percolation lattices.

The paper is organized as follows: In Section~\ref{perc_lat}, we define generalized percolation lattices and describe the cellular structure of such lattices in Theorem~\ref{lem_cel} and  in Section~\ref{pf1} we study outermost boundaries of star and plus connected components in generalized percolation lattices. Next, in Section~\ref{vac_cyc_sec}, we define the dual graph to a percolation lattice and describe conditions under which the dual graph has the properties of a percolation lattice. Following this, we use dual lattices to study vacant cycles of cells surrounding plus and star connected components. Finally in Section~\ref{pf_lr}, we study left right crossings of generalized percolation lattices in rectangles.

\setcounter{equation}{0}
\renewcommand\theequation{\arabic{section}.\arabic{equation}}
\section{Percolation lattices}\label{perc_lat}
Let~\(G = (V,E)\) be any connected finite planar graph in~\(\mathbb{R}^2\) where each edge is a straight line segment. Two vertices~\(v_1\) and~\(v_2\) are said to be adjacent if they share an edge in common. Two edges~\(e_1\) and~\(e_2\) are said to be adjacent if they share a vertex in common. A subgraph~\(P = (v_1,\ldots,v_l) \subseteq G\) is said to be a \emph{walk} if~\(v_i\) is adjacent to~\(v_{i+1}\) for~\(1 \leq i \leq l-1.\) If~\(e_i = (v_i,v_{i+1})\) is the edge containing end-vertices~\(v_i\) and~\(v_{i+1},\) we also represent~\(P = (e_1,\ldots,e_{l-1})\) and say that~\(v_1\) and~\(v_l\) are \emph{end-vertices} of~\(P.\) We say that~\(P\) is a \emph{circuit} if~\(P\) is a walk and~\(v_1 = v_l.\) We say that~\(P\) is a \emph{path} if~\(P\) is a walk and all the~\(l\) vertices in~\(P\) are distinct. Finally, we say that~\(P\) is a \emph{cycle} if~\(P\) is a path and~\(v_1 = v_l.\)

For a cycle~\(C \in G,\) let~\(A = A(C)\) be the bounded open set whose boundary is~\(C.\) We define the \emph{interior} of~\(C\) to be~\(A\) and the \emph{closed interior} of~\(C\) to be~\(A \cup C.\) We also define the \emph{exterior} of~\(C\) to be~\((C \cup A)^c\) and the \emph{closed exterior} of~\(C\) to be~\(A^c.\) We say that the graph~\(G\) is a \emph{percolation lattice} if every edge in~\(G\) belongs to a cycle. We say that a cycle~\(C\) in~\(G\) is a \emph{cell} if there exists no point of an edge of~\(G\) in the interior of~\(C.\) By definition any two distinct cells~\(C_1\) and~\(C_2\) have mutually disjoint interiors and the intersection~\(C_1 \cap C_2\) is either empty or a union of vertices and edges in~\(G.\) Any edge of~\(e\) belongs to at most two cells and we say that~\(e\) is \emph{unicellular} if there is at most one cell containing~\(e\) as an edge. 

\begin{figure}[tbp]
\centering
\includegraphics[width=2in, trim= 50 400 100 15, clip=true]{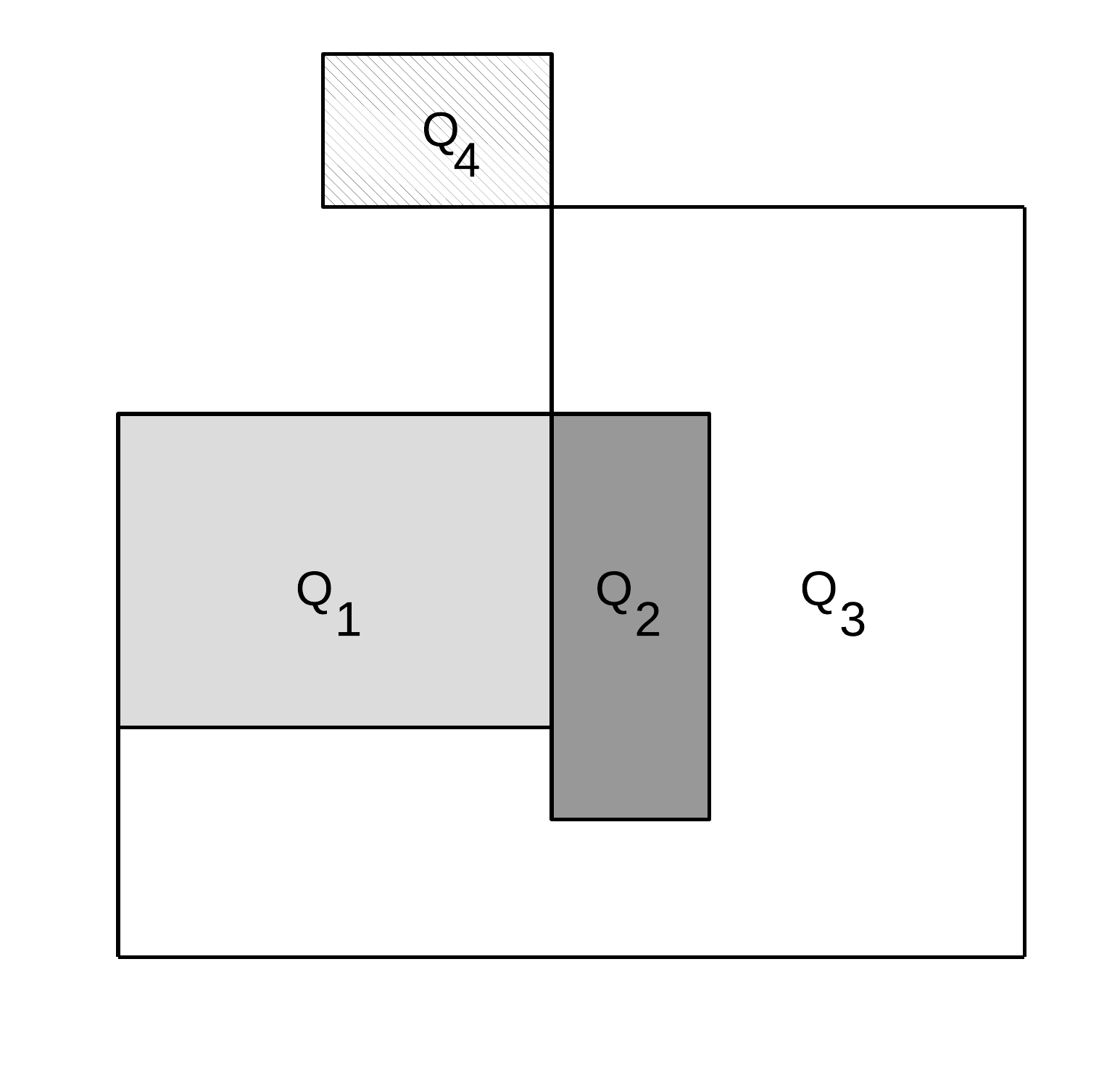}
\caption{Example of a percolation lattice~\(G\) along with the corresponding cellular decomposition~\(\bigcup_{i=1}^{4} Q_i.\) }
\label{perc_latt_fig}
\end{figure}

The following intuitive result captures the main features of percolation lattices as studied in~\cite{boll2}. 
\begin{Theorem}\label{lem_cel} If~\(G\) is a percolation lattice, then there are distinct cells\\\(Q_1,Q_2,\ldots,Q_T\) such that
\begin{equation}\label{cell_dec}
G = \bigcup_{i=1}^{T} Q_i.
\end{equation}
Moreover, the representation~(\ref{cell_dec}) is unique in the sense that if~\(V_1,\ldots,V_W\) are cells such that~\(G = \bigcup_{j=1}^{W} V_j,\) then~\(T = W\) and~\(\{V_i\}_{1 \leq i \leq T} = \{Q_i\}_{1 \leq i \leq T}.\)

The following additional properties hold:\\
\((x1)\) For every edge~\(e,\) there are at most two cells containing~\(e\) as an edge.\\
\((x2)\) If~\(e\) is contained in the closed interior of a cycle~\(C \in G,\) then there are two cells containing~\(e\) as an edge and both cells lie in the closed interior of~\(C.\) If~\(e\) is contained in the closed exterior of~\(C,\) then all cells containing~\(e\) lie in the closed exterior of~\(C.\)\\
\((x3)\) There are cycles~\(\Delta_1,\Delta_2,\ldots,\Delta_B\) with mutually disjoint interiors such that every cell in~\(G\) is contained in the closed interior of one of these cycles. For any~\(i \neq j,\) the cycles~\(\Delta_i\) and~\(\Delta_j\) have at most one vertex in common and an edge~\(e \in G\) is unicellular if and only if~\(e\) belongs to some cycle in~\(\{\Delta_l\}.\) 
\end{Theorem}
In Figure~\ref{perc_latt_fig}\((a),\) we illustrate the above result using a percolation lattice containing~\(4\) cells~\(Q_1,Q_2,Q_3\) and~\(Q_4.\)


\begin{figure}[ht!]
\centering
\begin{subfigure}{0.5\textwidth}
\centering
   \includegraphics[width=\textwidth, trim= 50 400 100 175, clip=true]{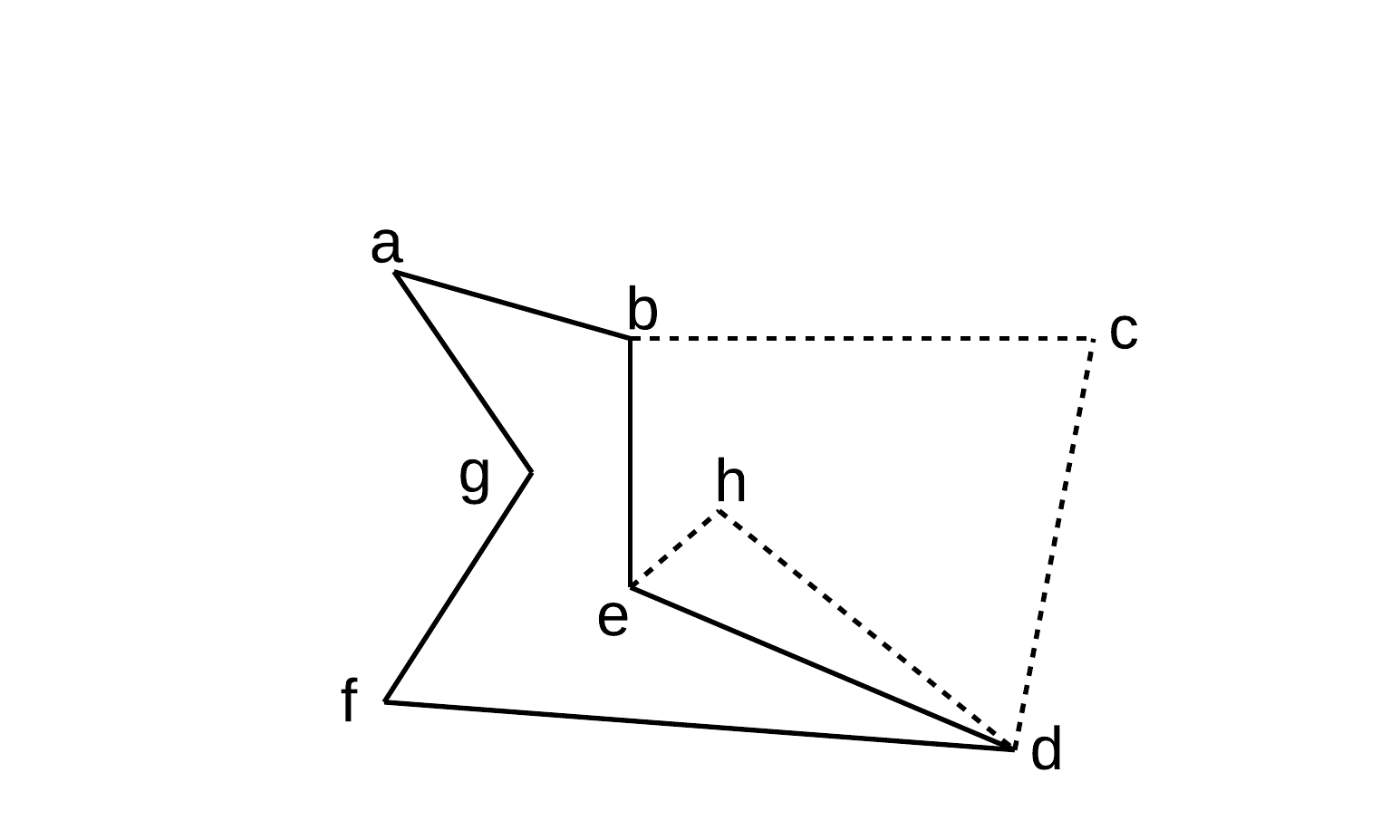}
  \caption{} 
\end{subfigure}
\begin{subfigure}{0.5\textwidth}
\centering
   \includegraphics[width=\textwidth, trim= 50 400 100 175, clip=true]{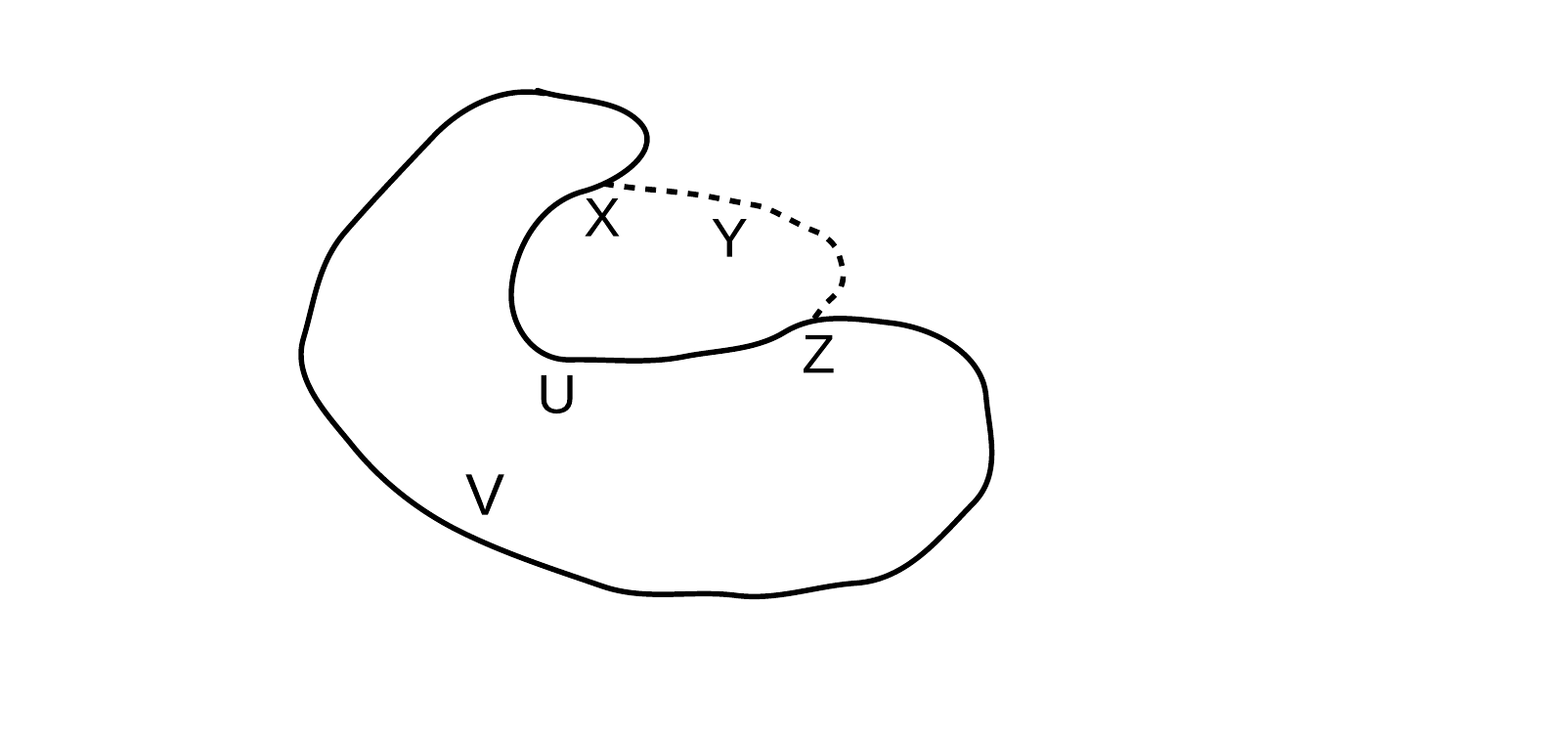}
   \caption{} 
  \end{subfigure}

\caption{\((a)\) The cycle~\(abcdfga\) formed by merging the cycles~\(abedfga\) and~\(behdcb.\) \((b)\) The path~\(P = XYZ \subset D\) lies in the exterior cycle~\(C = XUZVX.\)}
\label{cyc_fig}
\end{figure}

For completeness, we prove Theorem~\ref{lem_cel} using the following auxiliary result regarding merging of two cycles that is of independent interest and used throughout the paper.
\begin{Proposition}\label{thm1} Let \(C\) and \(D\) be cycles in the graph~\(G\) that have more than one vertex in common. There exists a unique cycle~\(E\) consisting only of edges of \(C\) and \(D\) with the following properties:\\
\((i)\) The closed interior of \(E\) contains the closed interior of both \(C\) and \(D.\)\\
\((ii)\) If an edge \(e\) belongs to \(C\) or \(D,\) then either \(e\) belongs to \(E\) or is contained in its closed interior.

Moreover, if \(D\) contains at least one edge in the closed exterior of \(C,\) then the cycle \(E\) also contains an edge of \(D\) that lies in the closed exterior of~\(C.\)
\end{Proposition}
The above result essentially says that if two cycles intersect at more than one point, there is an innermost cycle containing both of them in its interior.
For illustration, we refer to Figure~\ref{cyc_fig}\((a)\) where two cycles~\(abedfga\) and~\(behdcb\) have the edge~\(be\) and the vertex~\(d\) in common. The merged cycle~\(abcdfga\) contains both the smaller cycles in its closed interior. Analogous to~\cite{kest}, we use an iterative piecewise algorithmic construction for obtaining the merged cycle.


\emph{Proof of Proposition~\ref{thm1}}: Let~\(P \subset D\) be any path that has its end-vertices in the cycle~\(D_0 := C\) and lies in the exterior of~\(D_0\) (for illustration see Figure~\ref{cyc_fig}\((b)\) where~\(P= XYZ\) and~\(C = XUZVX\)). Letting~\(Q = XVZ \subset C\) the cycle~\(D_1 := P \cup Q\) then contains the cycle~\(C = D_0\) in its interior. We then repeat the above procedure with the cycle~\(D_1\) and look for another path~\(P_1 \subset D\) that lies in the exterior of~\(D_1\) and has end-vertices in~\(D_1.\) Arguing as before, we get a cycle~\(D_2\) that contains~\(P_1\) as a subpath and has the cycle~\(C\) in its closed interior. This procedure continues until we obtain a cycle~\(D_n\) that does not contain any edge of~\(D\) in its exterior.

We argue that~\(D_n\) is the desired cycle~\(E.\) By construction both~\(C\) and~\(D\) are contained in the closed interior of~\(D_n\) and~\(D_n\) consists of only edges of~\(C\) and~\(D\) so~\((i)\) and~\((ii)\) are true. Moreover, each cycle~\(D_i, 1 \leq i \leq n\) contains an edge of~\(D\) that lies in the closed exterior of~\(C.\) The uniqueness of~\(E\) is also true since if~\(E_1 \neq E\) is any edge satisfying~\((i)-(ii)\) and~\(E_1\) contains an edge of~\(E\) in its closed exterior, then some edge of~\(C\) or~\(D\) is present in the closed exterior of~\(E_1,\) a contradiction.~\(\qed\)



\emph{Proof of Theorem~\ref{lem_cel}}: Let~\(C_1,\ldots,C_T\) be the set of all cycles containing~\(e = (u,v)\) as an edge.  We shrink the cycles in an iterative manner as follows. Let~\(D_0 := C_1\) and suppose there exists some point of an edge of the cycle~\(C_j, j \geq 2\) in the \emph{interior} of~\(D_0.\) Because the graph~\(G\) is planar, there exists a path~\(P_1 \subset C_j\) present in the closed interior of~\(D_0.\) We shrink the cycle~\(D_0\) to the cycle~\(D_1\) containing the edge~\(e\) and the path~\(P_1.\) For illustration we again use Figure~\ref{cyc_fig}\((b)\) with~\(D_0 = C_1 = XYZVX\) and~\(C_2 = XYZVX.\) In this case~\(P_1 = XUZ\) and~\(D_0 = XYZVX.\) The ``shrinked" cycle~\(D_1 = XUZVX\) contains every edge~\(P_1.\) We now repeat the above procedure with the cycle~\(D_1\) and proceed iteratively to finally obtain a cycle~\(D_{fin}\) that does not contain any point of~\(\bigcup_{j=1}^{T} C_j\) in its interior.

The cycle~\(D_{fin} =: Q_e\) is a cell containing the edge~\(e.\) If there exists a cycle~\(C_j\) such that~\(C_j\) and~\(Q_e\) have mutually disjoint interiors, then we repeat the above procedure starting with the cycle~\(C_j\) and obtain another cell~\(R_e,\) also containing~\(e\) as an edge. By construction~\(Q_e\) and~\(R_e\) have mutually disjoint interiors and for any cycle~\(C_i\) we have that either~\(Q_e\) or~\(R_e\) is contained in the closed interior of~\(C_i.\) The set of all distinct cells in~\(\{Q_e\}_{e \in G} \bigcup \{R_e\}_{e \in G}\) is the desired cellular decomposition~(\ref{cell_dec}) of~\(G.\) The decomposition~(\ref{cell_dec}) is unique since every~\(V_i\) must necessarily be one of~\(Q_1,\ldots,Q_T.\)

The proof of~\((x1)\) and the proof of~\((x2)\) for~\(e\) present in the closed exterior of~\(C,\) follow from the above construction.  To prove the remaining part of~\((x2),\) suppose that~\(e\) is present in the closed interior of~\(C.\) Since~\(G\) is a percolation lattice, there exists a cycle~\(C_e \in G\) containing~\(e\) as an edge. This cycle~\(C_e\) contains an edge~\(e\) in the closed interior of~\(C\) and therefore a path~\(P_e \subset C_e\) with end-vertices~\(a,b \in C\) contained in the closed interior of~\(C.\)    Let~\(F_1\) and~\(F_2\) be the two paths with end-vertices~\(a,b\) that form the cycle~\(C.\) For reference, in Figure~\ref{cyc_fig}\((b),\) we have~\(X = a,Z=b, F_1 = XYZ, F_2 = XVZ\) and~\(P_e = XUZ.\) The cycles~\(P_e \cup F_1\) and~\(P_e \cup F_2\)
have mutually disjoint interiors and so arguing as before, we obtain two cells~\(Q_e\) and~\(R_e\) containing~\(e\) as an edge. By construction both~\(Q_e\) and~\(R_e\) are contained in the interior of~\(C.\)

The cycles in~\((x3)\) are obtained by repeatedly merging cells in~\(G\) as follows: We first pick one cell~\(Q_1\) and using Proposition~\ref{thm1}, merge~\(Q_1\) with another cell, say~\(Q_2,\) that shares an edge with~\(Q_1\) to get a new cycle~\(Q_{12}.\) We then pick another cell, say~\(Q_3,\) that shares an edge with~\(Q_{12}\) and lies in the exterior of~\(Q_{12}\) and merge these together to get a new cycle~\(Q_{123}.\) Continuing this way, we get a cycle~\(\Delta_1\) satisfying the property that no cell in~\(\{Q_j\}\) lying in the exterior of~\(\Delta_1\) shares an edge with~\(\Delta_1.\)
If there still exists cells in the exterior of~\(\Delta_1,\) then because~\(G\) is connected, one of these exterior cells (call it~\(Q_{21}\)) necessarily shares a vertex with~\(\Delta_1.\) We then repeat the above procedure starting with the cell~\(Q_{21}.\) Continuing this way until all cells are exhausted, we get the desired cycles~\(\Delta_l, 1 \leq l \leq B.\)

Finally, if~\(e\) is unicellular and is contained within~\(\Delta_l,\) then the cell~\(Q_e\) necessarily shares the edge~\(e\) with~\(\Delta_l.\) This completes the proof of~\((x3).\)~\(\qed\)

\setcounter{equation}{0}
\renewcommand\theequation{\arabic{section}.\arabic{equation}}
\section{Outermost boundaries}\label{pf1}
Let~\(G\) be a percolation lattice with cellular decomposition~\(G = \bigcup_{k=1}^{T} Q_i\) as in~(\ref{cell_dec}). We have the following definition of star and plus adjacency.
\begin{Definition}\label{nice_def} We say that two cells~\(Q_1\) and~\(Q_2\) in~\(G\) are \emph{star adjacent} if~\(Q_1 \cap Q_2\) contains a vertex in~\(G\) and \emph{plus adjacent} if~\(Q_1 \cap Q_2\) contains an edge in~\(G.\) 
\end{Definition}

We assign every cell~\(Q_k, 1 \leq k \leq T,\) one of the two states, occupied or vacant and assume that there exists an occupied cell~\(Q_0\) containing the origin. We say that the cell~\(Q_i\) is connected to the cell~\(Q_j\) by a \emph{star connected \(S-\)path} if there is a sequence of distinct cells~\((Y_1,Y_2,...,Y_t), Y_l \subset \{Q_k\}, 1 \leq l \leq t\) such that~\(Y_l\) is star adjacent to~\(Y_{l+1}\) for all \(1 \leq l \leq t-1\) and \(Y_1 = Q_i\) and \(Y_t = Q_j.\) If all the cells in \(\{Y_l\}_{1 \leq l \leq t}\) are occupied, we say that~\(Q_i\) is connected to~\(Q_j\) by an \emph{occupied} star connected \(S-\)path.

Let~\(C(0)\) be the collection of all occupied cells in~\(\{Q_k\}_{1 \leq k \leq T}\) each of which is connected to the cell~\(Q_0\) by an occupied star connected \(S-\)path. We say that \(C(0)\) is the  star connected occupied component containing the origin and let \(\{J_k\}_{1 \leq k \leq M} \subset \{Q_j\}\) be the set of all the occupied cells belonging to the component \(C(0).\)

\begin{figure}[ht!]
\centering
\begin{subfigure}{0.5\textwidth}
\centering
   \includegraphics[width=0.8\textwidth, trim= 50 350 100 175, clip=true]{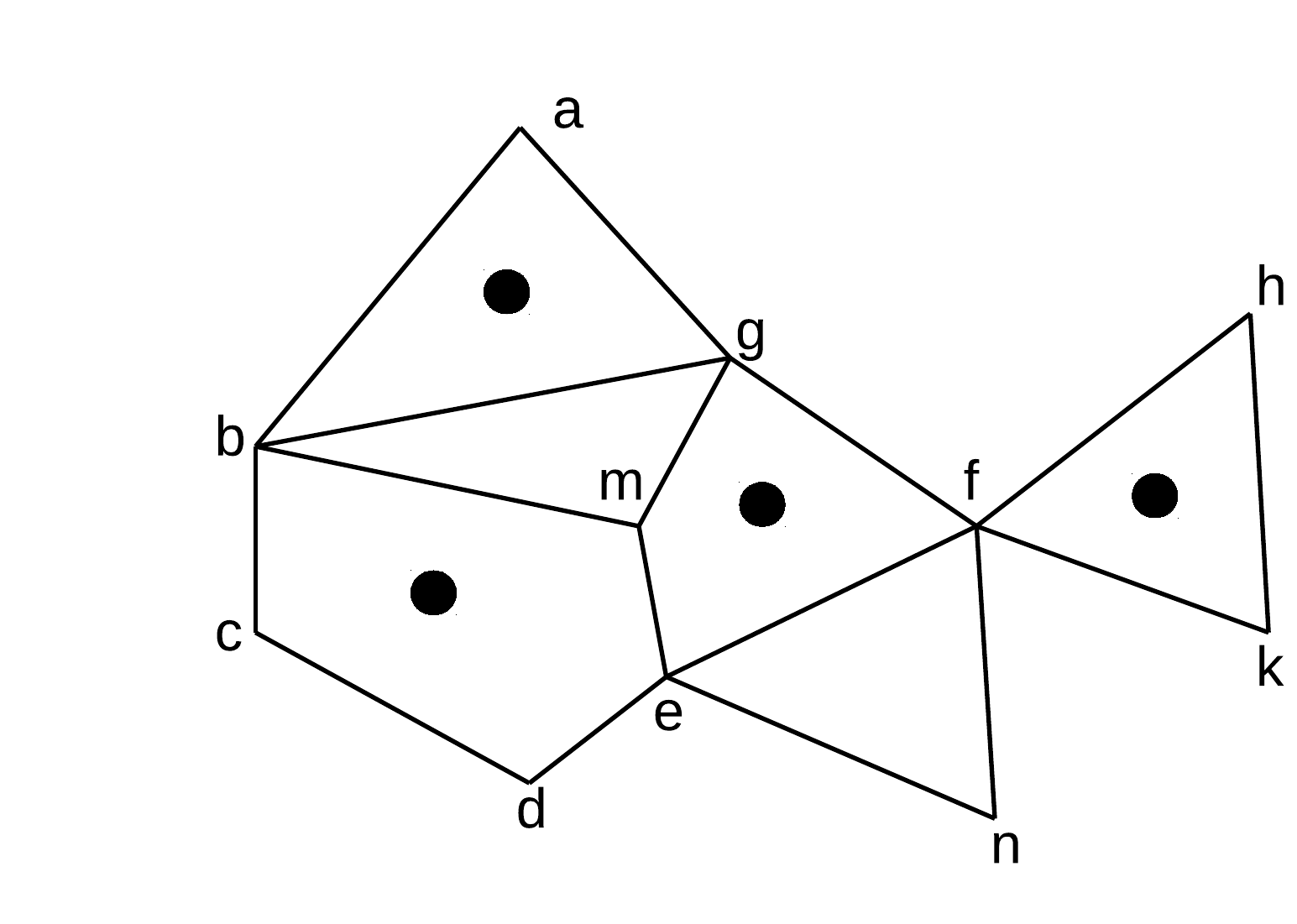}
  \caption{} 
\end{subfigure}
\begin{subfigure}{0.5\textwidth}
\centering
   \includegraphics[width=0.8\textwidth, trim= 20 300 150 150, clip=true]{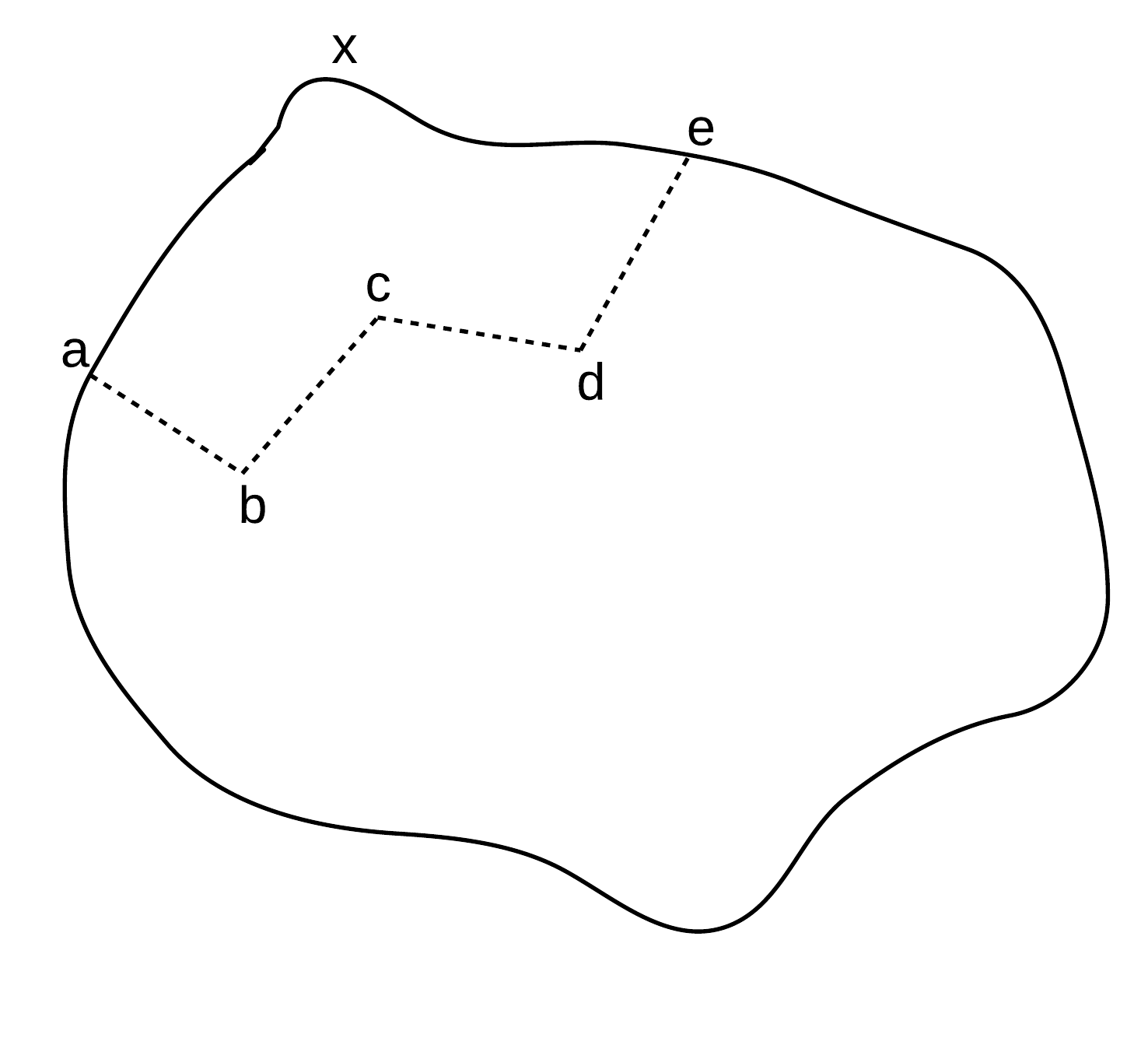}
   \caption{} 
  \end{subfigure}

\caption{\((a)\) A star connected component formed by the cells with circle inside them. The outermost boundary consists of two cycles~\(C_1 = abcdefga\) and~\(C_2 = fhkf.\) \((b)\) Replacing an interior path~\(P = abcde\) in the cycle~\(D_k\) (denoted by the wavy curve) with a path~\(Q  = axe \subseteq D_k.\)}
\label{fig_star_out}
\end{figure}

To define the outermost boundary of~\(C(0),\) we begin with a few preliminary definitions. Let \(G_0\) be the graph with vertex set and edge set, respectively being the vertex set and edge set of the cells~\(\{J_k\}_{1 \leq k \leq M} = C(0).\) An edge~\(e \in G_0\) is a  said to be \emph{boundary edge} if~\(e\) is unicellular or~\(e\) is adjacent to a vacant cell. (By definition,~\(e\) is already adjacent to an occupied cell of the component~\(C(0)\)). We have the following definition.
\begin{Definition} \label{out_def} We say that the edge \(e\) in the graph~\(G_0\) is an \emph{outermost boundary} edge of the component \(C(0)\) if the following holds true for every cycle \(C\) in~\(G_0:\) either \(e\) is an edge of~\(C\) or \(e\) is in the closed exterior of~\(C.\)

We define the outermost boundary \(\partial _0\) of \(C(0)\) to be the set of all outermost boundary edges of~\(G_0.\)
\end{Definition}
Thus outermost boundary edges cannot be contained in the interior of any cycle in the graph~\(G_0.\) We have the following result regarding the outermost boundary of the star component~\(C(0).\)
\begin{Theorem}\label{thm3} There are cycles~\(C_1,C_2,\ldots,C_n \subseteq G\) satisfying the following properties:\\
\((i)\) An edge~\(e \in \partial_0\) if and only if~\(e \in \bigcup_{i=1}^{n} C_i.\)\\
\((ii)\) The graph~\(\bigcup_{i=1}^{n} C_i\) is a connected subgraph of \(G_0.\)\\
\((iii)\) If \(i \neq j,\) the cycles~\(C_i\) and~\(C_j\) have disjoint interiors and have at most one vertex in common.\\
\((iv)\) Every occupied cell~\(J_k \in C(0)\) is contained in the interior of some cycle~\(C_{j}.\)\\
\((v)\) If~\(e \in C_{j}\) for some \(j,\) then \(e\) is a boundary edge belonging to an occupied square of~\(C(0)\) contained in the
interior of~\(C_j.\) If~\(e\) is not unicellular, then~\(e\) also belongs to a vacant cell lying in the exterior of all the cycles in~\(\partial_0.\)\\
Moreover, there exists a circuit \(C_{out}\) containing every edge of \(\cup_{1 \leq i \leq n} C_i.\)
\end{Theorem}
The outermost boundary \(\partial_0\) is a connected union of cycles satisfying properties~\((i)-(v)\) and is therefore an Eulerian graph with~\(C_{out}\) denoting the corresponding Eulerian circuit (see Chapter 1,~\cite{boll}). As an illustration, Figure~\ref{fig_star_out}\((a)\) describes a percolation lattice with six cells. The cells with circle inside them are occupied and the rest are vacant. The occupied cells form a star connected component and the outermost boundary consists of two cycles~\(C_1 = abcdefga\) and~\(C_2 = fhkf.\)

To prove Theorem~\ref{thm3}, we use the following Proposition also of independent interest.
\begin{Proposition}\label{outer} For every occupied cell~\(J_k \in C(0), 1 \leq k \leq M,\) there exists a unique cycle~\(D_k\) in~\(G_0\) satisfying the following properties.\\
\((a)\) The cell~\(J_k\) is contained in the closed interior of~\(D_k.\)\\
\((b)\) Every edge~\(e \in D_k\) is a boundary edge adjacent to an occupied cell of~\(C(0)\) present in the closed interior of~\(D_k.\) If~\(e\) is not unicellular, then~\(e\) is also adjacent to a vacant cell present in the closed exterior of~\(D_k.\)\\
\((c)\) If~\(C\) is any cycle in~\(G_0\) that contains~\(J_k\) in the interior, then every edge in~\(C\) either belongs to~\(D_k\) or is contained in the interior.
\end{Proposition}
Every edge of~\(D_k\) is also an outermost boundary edge in the graph~\(G_0\) and so we denote~\(D_k\) to be the \emph{outermost boundary cycle} containing the cell~\(J_k \in C(0).\) For example in Figure~\ref{fig_star_out}\((a),\) the outermost boundary cycle containing the cell~\(efgm\) is the cycle~\(C_1 = abcdefga.\)

Below we prove Proposition~\ref{outer} and Theorem~\ref{thm3} in that order.\\\\
\emph{Proof of Proposition~\ref{outer}}: Let \({\cal E} \neq \emptyset\) be the set of all cycles in the graph~\(G_0\) satisfying property~\((a);\) i.e., if \(C \in {\cal E}\) then~\(J_k\) in present in the closed interior of~\(C.\) The set~\({\cal E}\) is not empty since~\(J_k\) is itself a cycle containing~\(J_k\) in its closed interior and so belongs to~\({\cal E}.\) Moreover if~\(C_1\) and~\(C_2\) are any two cycles in~\({\cal E},\) then it cannot be the case that~\(C_1\) and~\(C_2\) have mutually disjoint interiors since both~\(C_1\) and~\(C_2\) contain the cell~\(J_k\) in its closed interior.
Thus it is possible to merge~\(C_1\) and~\(C_2\) using Proposition~\ref{thm1} to get a new cycle~\(C_3 \in {\cal E}\) that contains both~\(C_1\) and~\(C_2\)
in its closed interior. Continuing this way, we obtain an ``outermost cycle"~\(C_{fin}\) that contains all the cycles of~\({\cal E}\) in its closed interior.

By construction, the cycle~\(C_{fin}\)  satisfies properties~\((a)\) and~\((c).\) To see that~\((b)\) also holds,  suppose there exists an edge~\(e \in C_{fin}\) that is not a boundary edge. Since~\(e\) belongs to the graph~\(G_0,\) the edge~\(e\) is adjacent to an occupied cell~\(A_1 \in \{J_i\}.\) But since \(e\) is not a boundary edge, there exists one other cell~\(A_2 \in \{J_i\}\setminus A_1\) containing~\(e\) as an edge and moreover~\(A_2\) is also occupied. One of these cells, say \(A_1,\) is contained in the interior of~\(C_{fin}\) and the other cell~\(A_2,\) is contained in the exterior.

The cell~\(A_2\) and the cycle~\(C_{fin}\) have the edge \(e\) in common and thus more than one vertex in common. We then use Proposition~\ref{thm1} to obtain a larger cycle~\(C_{lar} \neq C_{fin}\) containing both~\(C_{fin}\) and~\(A_2\) in its closed interior. This is a contradiction to the fact that~\(C_{fin}\) satisfies property~\((c).\) Thus every edge~\(e\) of~\(C_{fin}\) is a boundary edge. By the same argument above, we also get that the edge~\(e\) cannot be adjacent to an occupied cell in the exterior of the cycle~\(C_{fin}.\) Thus~\(e\) is adjacent to an occupied cell in the interior of~\(C_{fin}\) and a vacant cell (if it exists) in the exterior.~\(\qed\)

\emph{Proof of Theorem~\ref{thm3}}: We argue that the set of distinct cycles in the set~\({\cal D} := \cup_{1 \leq k \leq M} \{D_k\}\) obtained in Proposition~\ref{outer} is the desired outermost boundary~\(\partial_0\) and satisfies the properties~\((i)-(v)\) mentioned in the statement of the theorem.

To prove~\((i),\) it suffices to see that every edge in the union of the cycles \({\cal D} = \cup_{1 \leq k \leq M} \{D_k\}\) is an outermost boundary edge. This is because, by definition, no edge with an end-vertex present in the interior of some cycle in~\({\cal D}\) can be an outermost boundary edge.
Now, suppose some edge~\(e \in D_k\) has an end-vertex in the interior of some cycle~\(C \subseteq G_0.\) This necessarily implies that at least one edge of~\(C\) is present in the exterior of~\(D_k\)  and moreover~\(C\) and~\(D_k\) cannot have a single common vertex. Therefore it is possible to merge~\(C\) and~\(D_k\) using Proposition~\ref{thm1} to get a bigger cycle~\(F_k \subset G\) containing~\(D_k\) in its closed interior, a contradiction to the construction of~\(D_k.\) This completes the proof of~\((i).\)

To prove~\((ii),\) we first see that the graph~\(G_0\) formed by the vertices and edges of the component~\(C(0),\) is connected. Indeed, let \(u_1\) and \(u_2\) be vertices in \(G_0\) so that each~\(u_i, i = 1,2\) is a corner of an occupied cell~\(J_i \in C(0).\) By definition, there is a star connected cell path connecting~\(J_1\) and~\(J_2,\) consisting only of cells in~\(C(0)\) and consequently, there exists a path in \(G_0\) from \(u_1\) to \(u_2.\) Now let~\(v_1\) and~\(v_2\) be vertices in~\({\cal D}\) belonging to cycles \(D_{r_1}\) and \(D_{r_2},\) respectively, for some \(1 \leq r_1,r_2 \leq M.\) By previous discussion, there exists a path~\(P_{12}\) from~\(v_1\) to~\(v_2\) containing only edges of~\(G_0\) and by construction, each such edge lies in the closed interior of some cycle in~\(\{D_k\}.\) This is true since all the occupied cells are present in some cycle in~\(\{D_k\}.\) For each sub-path~\(P \subseteq P_{12}\) that contains a point in the interior of some cycle~\(D_k\) and has end-vertices in~\(D_k,\) we replace~\(P\) with a path~\(Q \subseteq D_k\) (see Figure~\ref{fig_star_out}\((b)\)). Iteratively replacing all such interior paths, we get a path from~\(v_1\) to~\(v_2\) containing only edges of~\(\{D_k\}.\) Thus the union of cycles~\({\cal D}\) is connected and this proves~\((ii).\)

Property~\((iii)\) is true since otherwise we could merge the cycles~\(D_i\) and~\(D_j\) obtained in Proposition~\ref{outer} to get a larger cycle~\(D_{tot}\) containing both~\(D_i\) and~\(D_j\) in its closed interior, a contradiction to the fact that~\(D_i\) satisfies property~\((c)\) in Proposition~\ref{outer}. Indeed for any occupied cell~\(J_k \in C(0)\) the corresponding outermost boundary cycle~\(D_k\) satisfies property \((a)\) of Proposition~\ref{outer} and so~\((iv)\) is true. Moreover if~\(e \in D_k\) is any edge, then using the fact that~\(D_k\) satisfies property~\((b)\) of Proposition~\ref{outer}, we get that the edge \(e\) satisfies property~\((v).\)

Finally, to obtain the circuit~\(C_{out}\) containing the outermost boundary, we first compute the cycle graph \(T_{cyc}\) as follows. Let \(E_1,E_2,...,E_n\) be the distinct outermost boundary cycles in~\({\cal D}.\) Represent~\(E_i\) by a vertex~\(i\) in~\(T_{cyc}.\) If~\(E_i\) and~\(E_j\) share a corner, we draw an edge~\(e(i,j)\) between \(i\) and~\(j.\) Since the union of cycles~\(E_i, 1 \leq i \leq n\) is connected, we get that~\(T_{cyc}\) is connected as well.

Let~\(H_{cyc}\) be any spanning tree of~\(T_{cyc}\) and consider an increasing sequence of tree subgraphs~\(\{1\} = H_1 \subset H_2 \subset \ldots H_n = H_{cyc}.\) The graph \(H_1\) contains a single vertex \(\{1\}\) and so we set \(\Pi_1 = E_1\) to be the circuit obtained at the end of the first iteration. Having obtained the circuit~\(\Pi_i,\) let~\(q_{i+1} \in H_{i+1} \setminus H_i\) be adjacent to some leaf~\(v_{i} \in H_i.\) This implies that the cycle~\(E_{q_{i+1}}\) shares a vertex~\(w_i\) with the cycle~\(E_{v_i}.\) Since the circuit~\(\Pi_i\) contains~\(w_i,\) we assume that~\(w_i\) is the starting and ending vertex of~\(\Pi_i\) and also of~\(E_{q_{i+1}}.\) The concatenation of~\(\Pi_i\) and~\(E_{q_{i+1}}\) is the desired circuit~\(\Pi_{i+1}.\) Continuing this way, the final circuit~\(\Pi_{n}\) obtained is the desired circuit~\(C_{out}.\)~\(\qed\)

\subsection*{\em Plus connected components}
The techniques used in the previous sections also allows us to obtain the outermost boundary for plus connected components. We recall that cells~\(Q_i\) and~\(Q_j\) are said to be \emph{plus adjacent} if they share an edge between them. We say that the cell~\(Q_i\) is connected to the cell~\(Q_j\) by a \emph{plus connected \(S-\)path} if there is a sequence of distinct cells~\((Q_i = Y_1,Y_2,...,Y_t = Q_j) \subseteq \{Q_k\}_{1 \leq k \leq T}\) such that~\(Y_l\) is plus adjacent to~\(Y_{l+1}\) for all~\(1 \leq l \leq t-1.\) If all the cells in \(\{Y_l\}_{1 \leq l \leq t}\) are occupied, we say that~\(Q_i\) is connected to~\(Q_j\) by an \emph{occupied} plus connected \(S-\)path.

Let \(C^+(0)\) be the collection of all occupied cells in~\(\{Q_k\}_{1 \leq k \leq T}\) each of which is connected to the occupied cell~\(Q_0\) containing the origin, by an occupied plus connected~\(S-\)path. We say that \(C^+(0)\) is the plus connected occupied component containing the origin. Further we also define~\(G^+_0\) to be the graph with vertex set being the set of all vertices of the cells of~\(\{Q_k\}\) present in \(C^+(0)\) and edge set consisting of the edges of the cells of~\(\{Q_k\}\) present in~\(C^+(0).\)

Every plus connected component is also a star connected component and so the definition of outermost boundary edge in Definition~\ref{out_def} holds for the component~\(C^+(0)\) with~\(G_0\) replaced by~\(G^+_0.\) We have the following result.
\begin{Theorem}\label{thm2} The outermost boundary~\(\partial^+_0\) of~\(C^+(0)\) is unique cycle in~\(G^+_0\) with the following properties:\\
\((i)\) All cells of~\(C^+(0)\) are contained in the interior of~\(\partial^+_0.\)\\
\((ii)\) Every edge in~\(\partial^+_0\) is a boundary edge adjacent to an occupied cell of~\(C^+(0)\) contained in the interior of \(\partial^+_0\) and a vacant cell in the exterior.
\end{Theorem}
This is in contrast to star connected components which may contain multiple cycles in the outermost boundary.
In Figure~\ref{fig_star_out}\((a)\) for example, the union of the cells~\(bcdemb\) and~\(efgme\) forms a plus connected component whose outermost boundary is the cycle~\(bcdefgmb.\)

\emph{Proof of Theorem~\ref{thm2}}: Proposition~\ref{outer} holds with~\(C(0)\) replaced by~\(C^{+}(0).\) Since~\(C^+(0)\) is plus connected, the outermost boundary cycle~\(D_0\) the cell~\(Q_0\) in its interior also contains all the cells of~\(C^+(0)\) in its interior. Therefore cycle~\(D_0\) satisfies the conditions \((i)\) and \((ii)\) in the statement of the theorem, is unique and so~\(\partial_0^+ = D_0.\)~\(\qed\)

\renewcommand{\theequation}{\thesection.\arabic{equation}}
\setcounter{equation}{0}
\section{Vacant cell cycles surrounding occupied components} \label{vac_cyc_sec}
\begin{figure}[ht!]
\centering
\begin{subfigure}{0.5\textwidth}
\centering
\includegraphics[width=2.5in, trim = 20 330 120 140, clip = true]{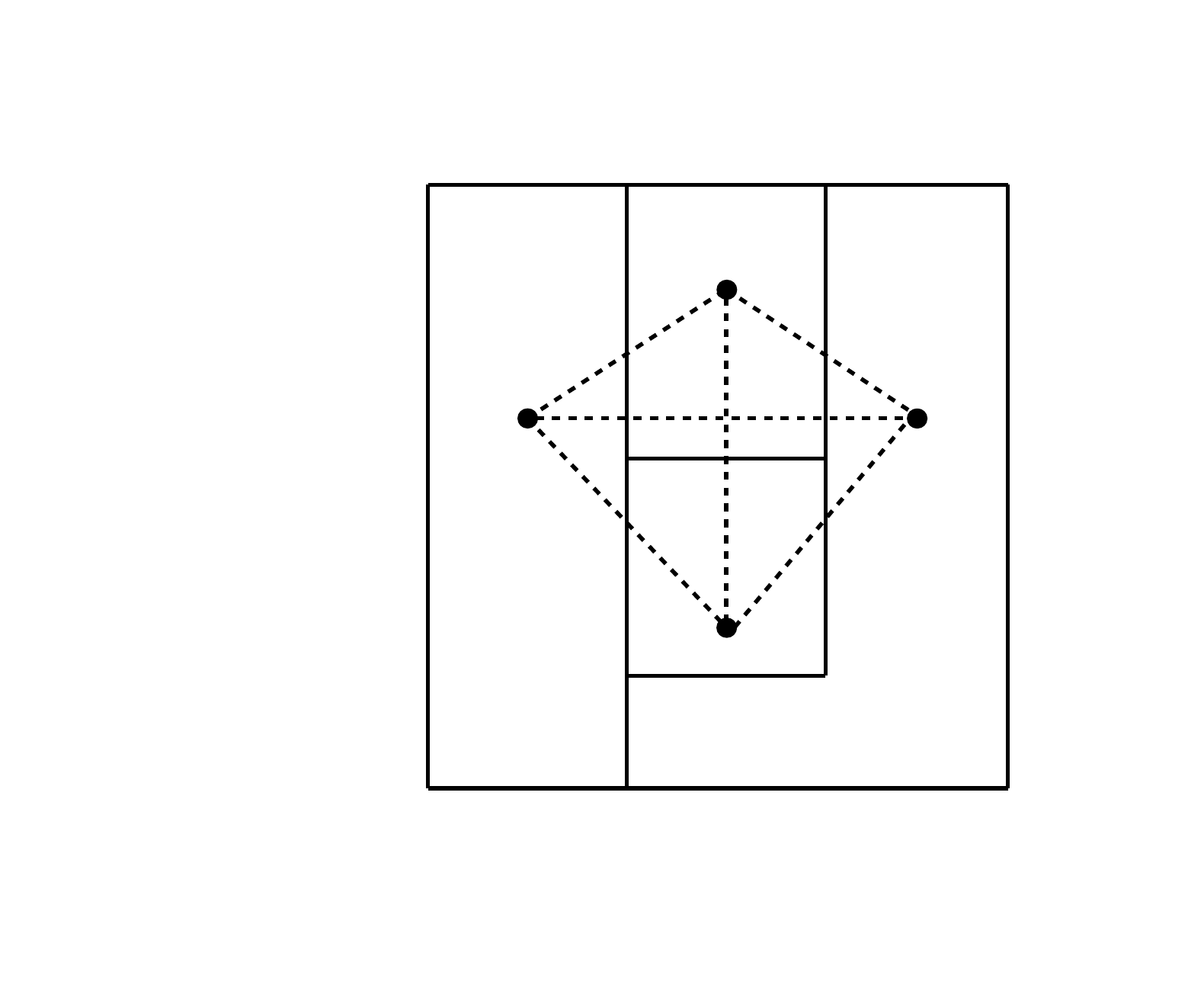}
  \caption{} 
\end{subfigure}
\begin{subfigure}{0.5\textwidth}
\centering
   \includegraphics[width=2.5in, trim= 20 330 120 140, clip=true]{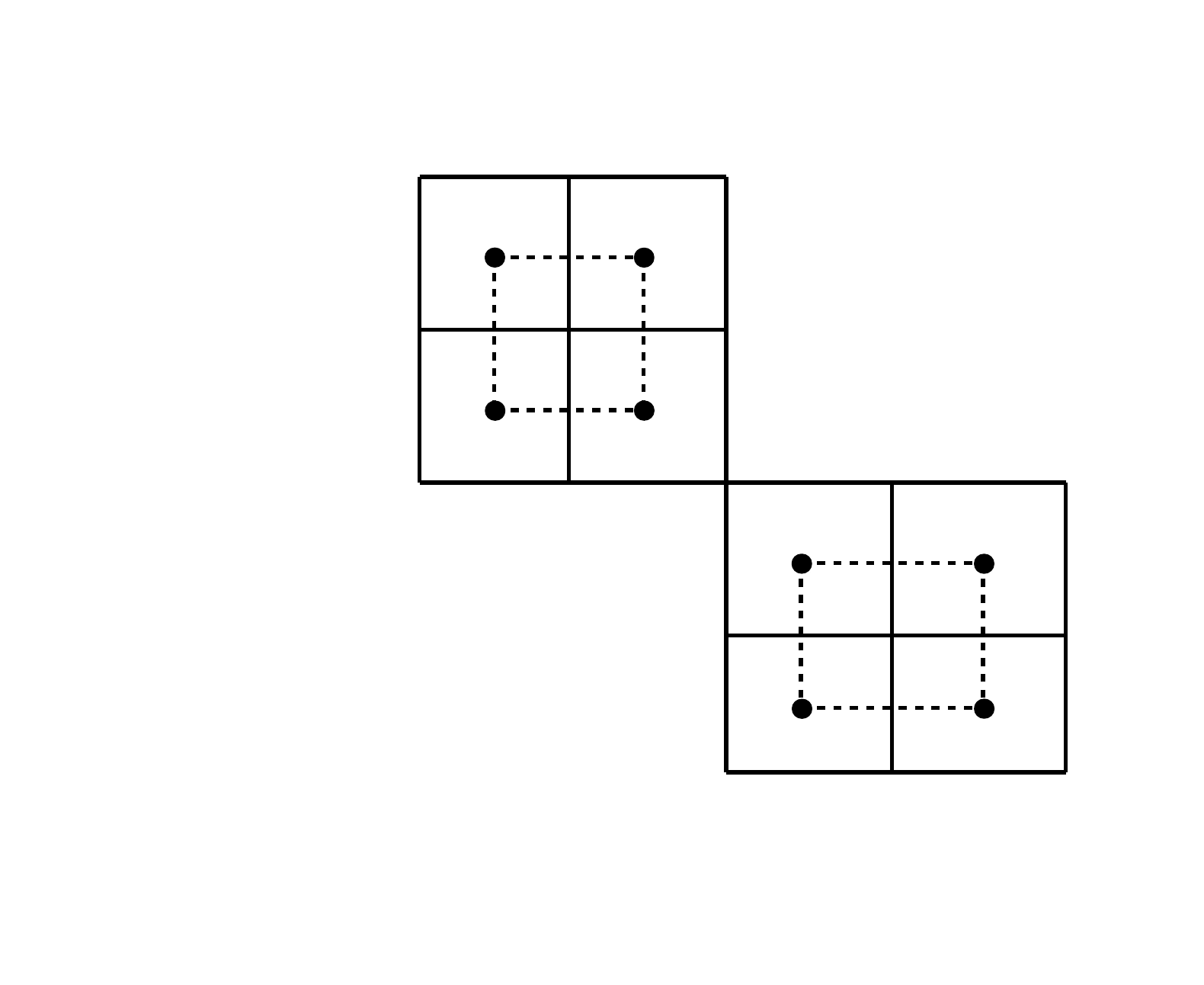}
   \caption{} 
  \end{subfigure}

\caption{\((a)\) An example where the dual graph~\(G_d\) is non-planar but connected and every edge of~\(G_d\) belongs to a cycle. \((b)\) An example where~\(G_d\) is planar and every edge of~\(G_d\) belongs to a cycle but~\(G_d\) is not connected.}
\label{dual_examples}
\end{figure}

In this section, we study vacant cycles of cells surrounding occupied star and plus components. We therefore begin with a discussion of the dual lattice. Let~\(G\) be any percolation lattice and let~\(G = \bigcup_{k=0}^{N} S_k\) be the cellular decomposition of~\(G.\)
\begin{Definition}\label{dual_def}
We say that a graph~\(G_d\) is \emph{dual} to~\(G\) if the following conditions hold:\\
\((d1)\) Every cell in~\(G\) contains exactly one vertex of~\(G_d\) in its interior.\\
Suppose vertex~\(w \in G_d\) is the present in the interior of the cell~\(Q(w)\) of~\(G.\)\\
\((d2)\) Vertices~\(w_1,w_2 \in G_d\) are adjacent if and only if the cells~\(Q(w_1)\) and~\(Q(w_2)\) are plus adjacent.
\end{Definition}
To study the similarity between~\(G\) and~\(G_d,\) we would like to first ensure that the dual graph~\(G_d\) is also a percolation lattice, i.e., we prefer that~\(G_d\) is connected, planar and each edge of~\(G_d\) belongs to a cycle. This is because, there are in fact dual graphs that satisfy exactly two of these three properties. For example, consider two plus adjacent squares~\(S_1\) and~\(S_2\) of same side length, to be the graph~\(G.\) If we let the centres of~\(S_1\) and~\(S_2\) be the vertex set of~\(G_d,\) then the edge set of~\(G_d\) is the single edge joining the centres of~\(S_1\) and~\(S_2.\) The graph~\(G_d\) is connected and planar but acyclic.
In Figure~\ref{dual_examples}~\((a),\) we have an example of a graph~\(G\) (denoted by solid lines) and the corresponding dual graph~\(G_d\) (denoted by the dotted lines). The graph~\(G_d\) is connected and every edge of~\(G_d\) belongs to a cycle but~\(G_d\) is not planar. In Figure~\ref{dual_examples}~\((b),\) the dual graph~\(G_d\) is planar and every edge of~\(G_d\) belongs to a cycle but~\(G_d\) is not connected.

Throughout the paper, we assume that the graph~\(G\) admits a dual graph~\(G_d\) satisfying the following properties:\\
\((a1)\) (Niceness property) The percolation lattice~\(G\) is nice in the sense that any two plus adjacent cells in~\(G\) share exactly one edge in common and no other vertex.\\
\((a2)\) (Interior edge property) Any edge~\((w_1,w_2) \in G_d\) is present in the interior of the cycle formed by merging the plus adjacent cells~\(Q(w_1)\) and~\(Q(w_2).\)\\
\((a3)\) The dual graph~\(G_d\) is a connected and planar percolation lattice and~\(G\) is dual to~\(G_d.\)\\
\((a4)\) The dual graph~\(G_d\) satisfies the niceness and interior edge property.\\
\((a5)\) If~\(z\) is a vertex of the dual cell~\(R(v)\) and~\(Q(z)\) is the cell in~\(G\) containing~\(z,\) then~\(v\) is a vertex of~\(Q(z).\)\\
Thus~\(G\) and~\(G_d\) have similar cellular structure. For example, the usual square lattice on the plane satisfies properties~\((a1)-(a5).\)

Let~\(G \) be a percolation lattice with cellular decomposition~\(\bigcup_{k=0}^{N} S_k\) and let~\(G_d\) be a lattice dual to~\(G\) satisfying properties~\((a1)-(a5)\) and having cellular decomposition~\(G_d = \bigcup_{l=0}^{L}W_l.\) We say that the sequence~\(L_S = (S_1,...,S_m)\) is a \emph{plus connected cell path} in~\(G\) if for each \(1 \leq i \leq m-1,\) the cell~\(S_i\) is plus adjacent with the cell~\(S_{i+1}.\) We say that \(L_S\) is a \emph{plus connected cell cycle} if~\(S_i\) is adjacent to~\(S_{i-1}\) and~\(S_{i+1}\) for each~\(1 \leq i \leq m\) with the notation that~\(S_{m+1} = S_1.\) Analogous definition holds for star connected paths and cycles.

\begin{figure}[tbp]
\centering
\includegraphics[width=2.5in, trim = 20 320 120 260, clip = true]{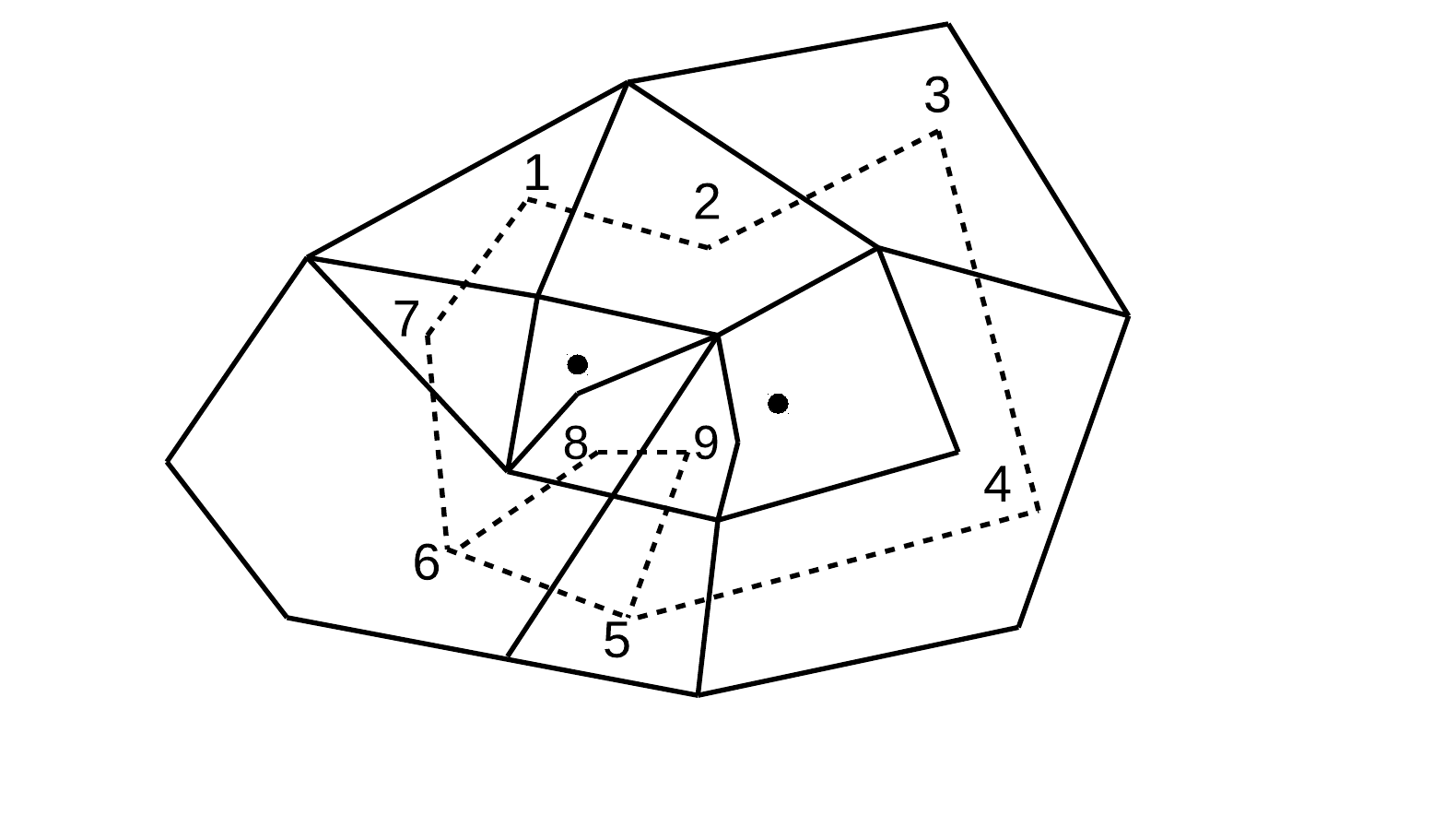}
\caption{Cycle of vacant cells surrounding a star connected occupied component.}
\label{vac_cyc_fig}
\end{figure}

Let~\(C(0) = \bigcup_{k=0}^{T} J_k \) be the star connected occupied cell component containing the cell~\(S_0\) with origin in its interior.
By definition, every cell in~\(C(0)\) is connected to~\(S_0\) by a star connected cell path.
We have the following result.
\begin{Theorem}\label{thm4} Suppose properties~\((a1)-(a5)\) hold and suppose every vertex in the component~\(C(0)\) is present in the interior of some dual cell of~\(G_d.\)

There exists a unique cycle~\({\cal P}_{out} = (w_1,\ldots,w_s) \subset G_d\) such that each vertex~\(w_i\) is present in the interior of a vacant cell~\(V_i\) and satisfies the following properties:\\
\((i)\) For every~\(i, 1 \leq i \leq s,\) the cell~\(V_i\) is vacant and star adjacent to some occupied cell in~\(C(0).\)\\
\((ii)\) All occupied cells in~\(C(0)\) are contained in the interior of~\({\cal P}_{out}.\)\\
\((iii)\) If~\(F_{out} \neq {\cal P}_{out}\) is any  other cycle in~\(G_d\) that satisfies \((i)-(ii)\) above, then~\(F_{out}\) is contained in the closed interior of~\({\cal P}_{out}.\)
\end{Theorem}
The sequence of cells in~\((V_1,\ldots,V_s)\) form a plus connected cycle of vacant cells surrounding the star connected component~\(C(0).\)
In Figure~\ref{vac_cyc_fig}, the two cells containing the circles form the star connected occupied component. Every other cell is vacant. The two dual cycles~\(12345671\) and~\(1234598671\) both satisfy~\((i)-(ii)\) and~\({\cal P}_{out} = 12345671.\)

\emph{Proof of Theorem~\ref{thm4}}: Let \(\partial_0\) denote the outermost boundary of the star connected component~\(C(0)\) in the percolation lattice~\(G.\) From Theorem~\ref{thm3} we have that~\(\partial_0 = \cup_{1 \leq i \leq n}C_i\) is a connected union of cycles \(\{C_i\}\) with mutually disjoint interiors and moreover, for~\(i \neq j,\) the cycles~\(C_i\) and~\(C_j\) have at most one common vertex.

If vertices~\(v_1, v_2 \in G\) are adjacent in the outermost boundary~\(\partial_0\) of the component~\(C(0),\) then the corresponding dual cells~\(R(v_1)\) and~\(R(v_2)\) are plus adjacent (property~\((a4)\). Also, because~\(\partial_0\) is connected (see property~\((ii)\) Theorem~\ref{thm3}), the union of dual cells
\begin{equation}\label{cv_par_def}
C_V(\partial_0) := \bigcup_{v \in \partial_0} R(v)
\end{equation}
is a plus connected component in the dual graph~\(G_d.\) Moreover, each edge of~\(\partial_0\) is present in the interior of closed interior of the union of cells of~\(C_V(\partial_0).\) Using Theorem~\ref{thm2} we therefore have that the outermost boundary~\(\partial_V(\partial_0)\) of~\(C_V(\partial_0)\) is a single cycle in~\(G_d\) containing all cells of~\(C_V(\partial_0)\) in its closed interior and all edges of~\(\partial_0\) in its interior. Here the dual outermost boundary~\(\partial_V(\partial_0)\) is obtained as follows. Every dual cell belonging to~\(C_V(\partial_0)\) is labelled~\(1\) and every dual cell sharing an edge with a cell in~\(C_V(\partial_0)\) and not belonging to~\(C_V(\partial_0),\) is labelled~\(0.\) We then apply Theorem~\ref{thm2} with label~\(1\)  cells as occupied and label~\(0\) cells as vacant.

Suppose \(z_1,z_2,\ldots,z_t\) are the vertices of the dual cycle \(\partial_V(\partial_0)\) encountered in that order; i.e., the vertex~\(z_1\) is adjacent to~\(z_2,\) the vertex~\(z_2\) is adjacent to~\(z_3\) and so on. Each vertex~\(z_i\) is a vertex of the dual cell~\(R(v)\) for some~\(v \in \partial_0.\) Therefore if~\(Q(z_i)\) is the cell in~\(G\) containing~\(z_i\) in its interior, then~\(v\) is a vertex of~\(Q(z_i),\) by property~\((a5).\) Moreover~\(Q(z_i)\) lies in the exterior of~\(\partial_0\) and is adjacent to the vertex~\(v\) of~\(C(0).\) This implies that~\(Q(z_i)\) must necessarily be vacant. This proves that the cycle~\(\partial_V(\partial_0)\) satisfies properties~\((i)-(ii).\)

To get a unique dual cycle satisfying properties~\((i)-(iii),\) we merge all dual cycles satisfying properties~\((i)-(ii).\) This is possible since any two dual cycles satisfying~\((i)-(ii)\) both contain the cell~\(S_0\) in their respective interiors and therefore cannot have mutually disjoint interiors.~\(\qed\)


\subsection*{\em Plus connected components}
In this subsection we let~\(C^+(0) = \bigcup_{k=0}^{T} J_k \) be the plus connected occupied cell component containing the cell~\(S_0\) with origin in its interior. By definition, every cell in~\(C^+(0)\) is connected to~\(S_0\) by a plus connected cell path and the outermost boundary~\(\partial_0^+\) of~\(C^+(0)\) is a single cycle containing all the cells of~\(C^+(0)\) in its interior.
We have the following result.
\begin{Theorem}\label{thm55} There exists a star connected cell cycle~\({\cal M}_{out} = (Y_1,\ldots,Y_t) \subset G\) such that:\\
\((i)\) Each cell~\(Y_i\) is vacant, lies in the exterior of~\(\partial_0^+\) and is plus adjacent to some occupied cell of~\(C^+(0).\)\\
\((ii)\) The outermost boundary of~\({\cal M}_{out}\) is a single cycle containing all the cells of~\({\cal M}_{out} \cup C^+(0)\) in its interior.
\end{Theorem}
The sequence of cells in~\((Y_1,\ldots,Y_t)\) form a star connected cell cycle of vacant cells surrounding the plus connected component~\(C^+(0).\)

\emph{Proof of Theorem~\ref{thm55}}: From Theorem~\ref{thm2}, we have that the outermost boundary~\(\partial^+_0 = (e_1,\ldots,e_t)\) of~\(C^+(0)\) is a single cycle containing all cells of~\(C^+(0)\) in its interior. Moreover every edge~\(e_i\) of~\(\partial^+(0)\) belongs to a occupied cell of~\(C^+(0)\) and also to a vacant cell~\(Z_i\) lying in the exterior of~\(\partial^+_0.\) It is possible that multiple edges in~\(\partial_0^+\) belong to the same cell~\(Z_i\) and so the sequence~\((Z_1,Z_2,\ldots,Z_m)\) could have repetitions. In such a case, we remove recurring entries and assume without loss of generality that~\(Z_i\) is star adjacent to~\(Z_{i+1}\) for~\(1 \leq i \leq m\) with the notation that~\(Z_{m+1} = Z_1.\)

The set of cells in~\(\Gamma := (Z_1,\ldots,Z_m)\) form a star connected component and we suppose that the sequence of cells~\((Z_1,\ldots,Z_{n})\) form a star connected cycle~\(L_Z\) for some~\(n \leq m.\) The outermost boundary~\(\partial_Z\) of~\(L_Z\) is a connected union of cycles~\(\bigcup_{1 \leq i \leq T} D_i\) and contains every cell of~\(L_Z\) in the interior of some~\(D_i.\) If some cell of~\(C^+(0)\) is contained in the interior of~\(D_i,\) then because~\(C^+(0)\) is plus connected, every cell of~\(C^+(0)\) is contained in the interior of~\(D_i.\) Moreover, since~\(D_i\) and~\(D_j\) share at most one vertex in common, it must the case that~\(T = 1\) and so the outermost boundary of the star cycle~\((Z_1,\ldots,Z_n)\) is the single cycle~\(D_1.\)

If on the other hand, every cell of~\(C^+(0)\) is contained in the exterior of every cycle of~\(\partial_Z,\) then \emph{every} edge of~\(\partial_Z\) is the edge of some occupied cell of~\(C^+(0)\) that lies in the exterior of~\(\partial_Z.\) This implies that the outermost boundary~\(\partial_0^+\) of~\(C^+(0)\) is contained in the \emph{strict} exterior of every cycle in~\(\partial_Z.\) But this contradicts the fact that each~\(Z_i\) contains at least one edge of~\(\partial_0^+\) and so there is at least one edge of~\(\partial_0^+\) contained in the closed interior of some cycle of~\(\partial_Z.\)~\(\qed\)

\renewcommand{\theequation}{\thesection.\arabic{equation}}
\setcounter{equation}{0}
\section{Left right and top bottom crossings in rectangles}\label{pf_lr}
In this section, we study the mutual exclusivity of left right and top down crossings in a rectangle. As before we assume that the percolation lattice~\(G = \bigcup_{k=0}^{N}S_k\) and the dual lattice~\(G_d = \bigcup_{l=0}^{M}W_l\) satisfy properties~\((a1)-(a5)\) and the origin is the present in the interior of the cell~\(S_0.\)

\begin{figure}[tbp]
\centering
\includegraphics[width=3in, trim = 40 220 20 170, clip = true]{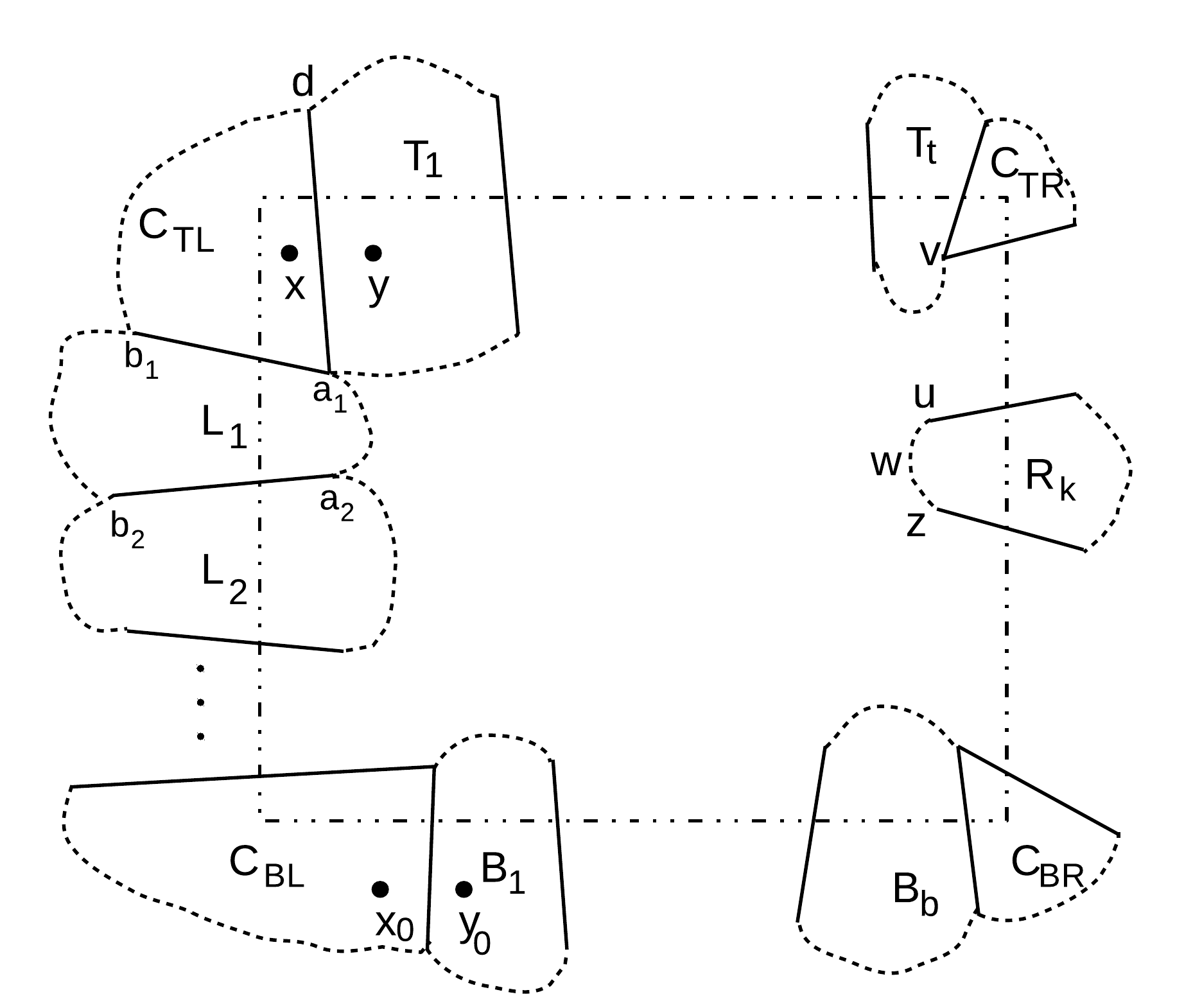}
\caption{The sides of the rectangle are nicely covered by the cells of~\(G.\)}
\label{nice_cover}
\end{figure}

For a fixed rectangle~\(R,\) we assume that the sides of~\(R\) are \emph{nicely covered} by cells of~\(G\) as shown in Figure~\ref{nice_cover}. Consider the edges~\(a_1b_1\) and~\(a_2b_2\) intersecting the left side of~\(R.\) The vertices~\(a_1\) and~\(a_2\) are connected by a path~\(A_1\) (shown by dotted line) in the interior of~\(R.\) Similarly the vertices~\(b_1\) and~\(b_2\) are connected by a path~\(B_1\) in the exterior of~\(R.\)  The union~\(A_1 \cup B_1 \cup \{a_1b_1,a_2b_2\}\) forms the cell~\(L_1.\) We define~\(L_1,\ldots,L_l\) to be the left cells. Similarly, the cells~\(T_1,\ldots,T_t\) are top cells, the right cells are~\(R_1,\ldots,R_r\) and the bottom cells are~\(B_1,\ldots,B_b.\) Any cell contained in the closed interior of~\(R\) is called an \emph{interior} cell. 

The cells~\(C_{TL},C_{TR}, C_{BL}\) and~\(C_{BR}\) each contain a corner of the rectangle and have a single vertex contained in the interior of the rectangle. These four cells, called \emph{corner cells}, are not plus adjacent to any interior cell. As in Figure~\ref{nice_cover}, we assume that the cell~\(L_1\) is plus adjacent to~\(L_2\) and~\(C_{TL}\) and does not share a vertex with any other left, right,top or bottom cell. We make analogous assumptions for each of the left, right, top and bottom cells.

Finally, we also assume that the cells are \emph{nicely padded} in the following way: Let~\({\cal L}_{top}\) be the infinite line containing the top side of~\(R\) and define~\({\cal L}_{bottom}, {\cal L}_{left}\) and~\({\cal L}_{right}\) analogously. If~\(Z_1\) is any cell intersecting~\({\cal L}_{top}\) and star adjacent to a cell intersecting~\(R\) and~\(Z_2\) is any cell intersecting~\({\cal L}_{bottom}\) and star adjacent to a cell intersecting~\(R,\) then~\(Z_1\) and~\(Z_2\) are not star adjacent. A similar assumption holds for cells intersecting~\({\cal L}_{left}\) and~\({\cal L}_{right}.\)

Assuming that~\(R\) is nicely covered and padded as described above, we have the following definition of left right crossings.
\begin{Definition}\label{lr_crossing_def}
A plus connected cell path~\(L = (J_1,\ldots,J_m) \subset \{S_k\}, m \geq 3\) is said to be a \emph{plus connected left right crossing} of a rectangle~\(R\) if~\(J_1\) is a left cell, the cell~\(J_m\) is a right cell and every~\(J_i, 2 \leq i \leq m-1\) is an interior cell.
\end{Definition}
By the nicely padded assumption,~\(L\) must contain at least one interior cell. Every interior cell is now assigned one of the following two states: occupied or vacant. If every interior cell in a left right crossing~\(L\) of~\(R\) is occupied, we say that \(L\) is an \emph{occupied plus connected left right crossing} of the rectangle \(R.\) We denote \(LR^+(R,O)\) and \(LR^+(R,V)\) to be the events that the rectangle~\(R\) contains an occupied and vacant plus connected left right crossing, respectively.

We have a similar definition of plus connected top down crossings and denote~\(TD^+(R,O)\) and \(TD^+(R,V)\) to be the events that the rectangle~\(R\) contains an occupied and vacant plus connected top down crossing, respectively. Replacing plus adjacent with star adjacent, we obtain an analogous definition for star connected left right and top down crossings. We have the following result.
\begin{Theorem}\label{thm_lr} Suppose properties~\((a1)-(a5)\)  hold and~\(R\) is nicely covered and nicely padded by cells as in Figure~\ref{nice_cover}. Also suppose that every vertex of a cell intersecting~\(R\) is present in the interior of a dual cell and that there is at least one plus connected left right crossing and one plus connected top bottom crossing. We have the following.\\
\((i)\) One of the events \(LR^+(R,O)\) or \(TD^*(R,V)\) always occurs but not both.\\
\((ii)\) One of the events \(LR^*(R,O)\) or \(TD^+(R,V)\) always occurs but not both.
\end{Theorem}
The above result describes the mutual exclusivity of occupied and vacant left right and top down crossings in any rectangle.

\subsection*{\em Proof of Theorem~\ref{thm_lr}}
We prove the following three statements.\\
\((I)\) The events~\(LR^*(R,O)\) and~\(TD^+(R,V)\) cannot both occur simultaneously.\\
\((II)\) If~\(LR^*(R,O)\) does not occur, then~\(TD^+(R,V)\) must necessarily occur.\\
\((III)\) If~\(LR^+(R,O)\) does not occur, then~\(TD^*(R,V)\) occurs.\\
Using~\((I)-(III)\) and the fact that a top down crossing of~\(R\) is a left right crossing of the rectangle~\(R'\) obtained by rotating~\(R\) by ninety degrees around the centre, we then get Theorem~\ref{thm_lr}.\\

\emph{Proof of~\((I)\)}: Suppose there exists a star connected occupied left right crossing~\(L_s\) of~\(R\) and let~\(\Gamma = (e_1,\ldots,e_t)\) be a path in~\(L_s\) crossing~\(R\) from left to right so that~\(e_1\) intersects the left edge of~\(R,\) the edge~\(e_t\) intersects the right edge of~\(R\) and each~\(e_i, 2 \leq i \leq t-1\) belongs to an interior occupied cell of~\(L_s.\) The path~\(\Gamma\) divides the rectangles into two halves.

\begin{figure}[ht!]
    \centering
    \begin{subfigure}[t]{0.4\textwidth}
    \centering  
        \includegraphics[width=\textwidth, trim = 40 320 20 150, clip = true]{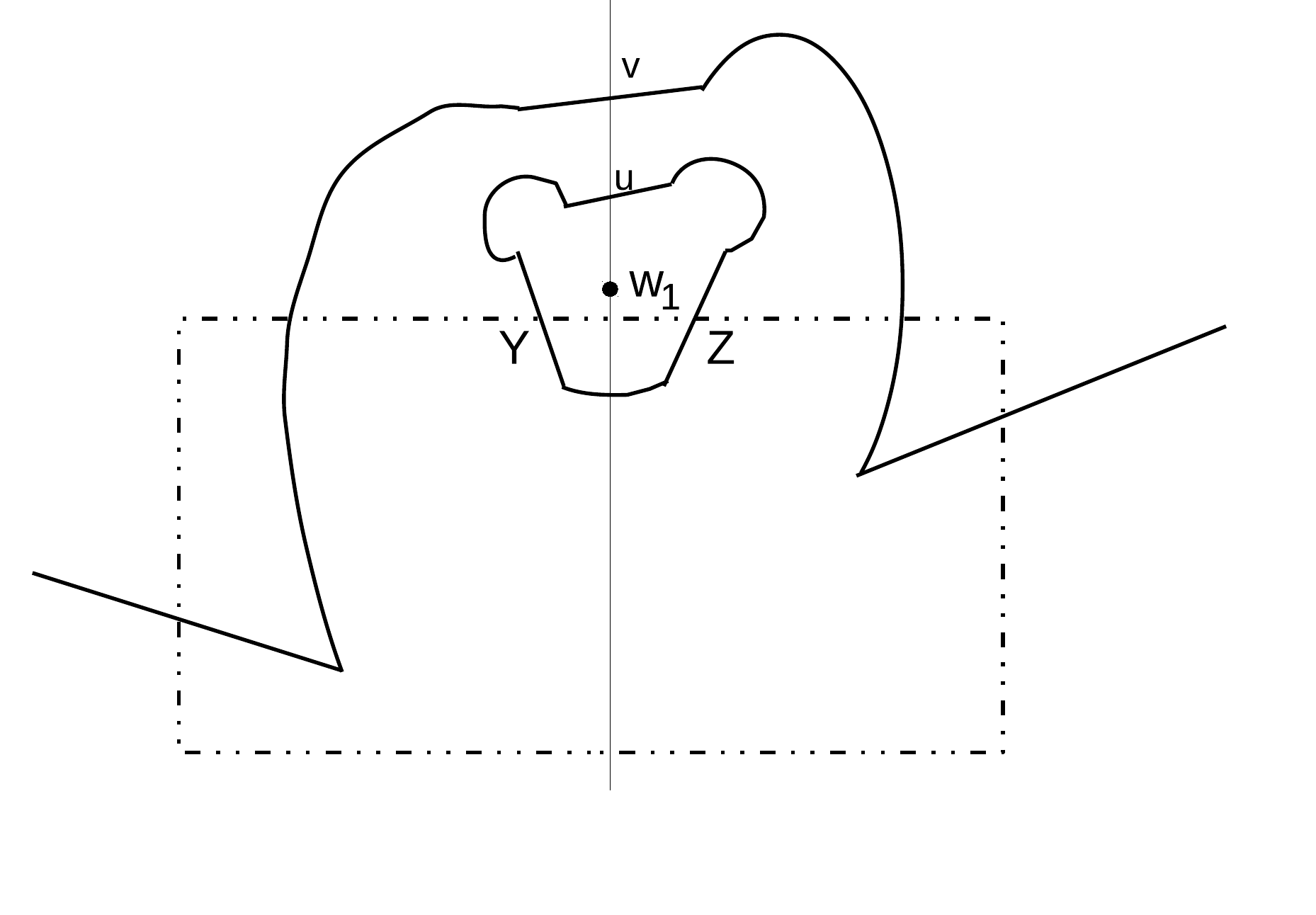}
        \caption{}
     \end{subfigure} 
     ~
     \begin{subfigure}[t]{0.4\textwidth}
     \centering
        \includegraphics[width=\textwidth, trim = 40 320 20 150, clip = true]{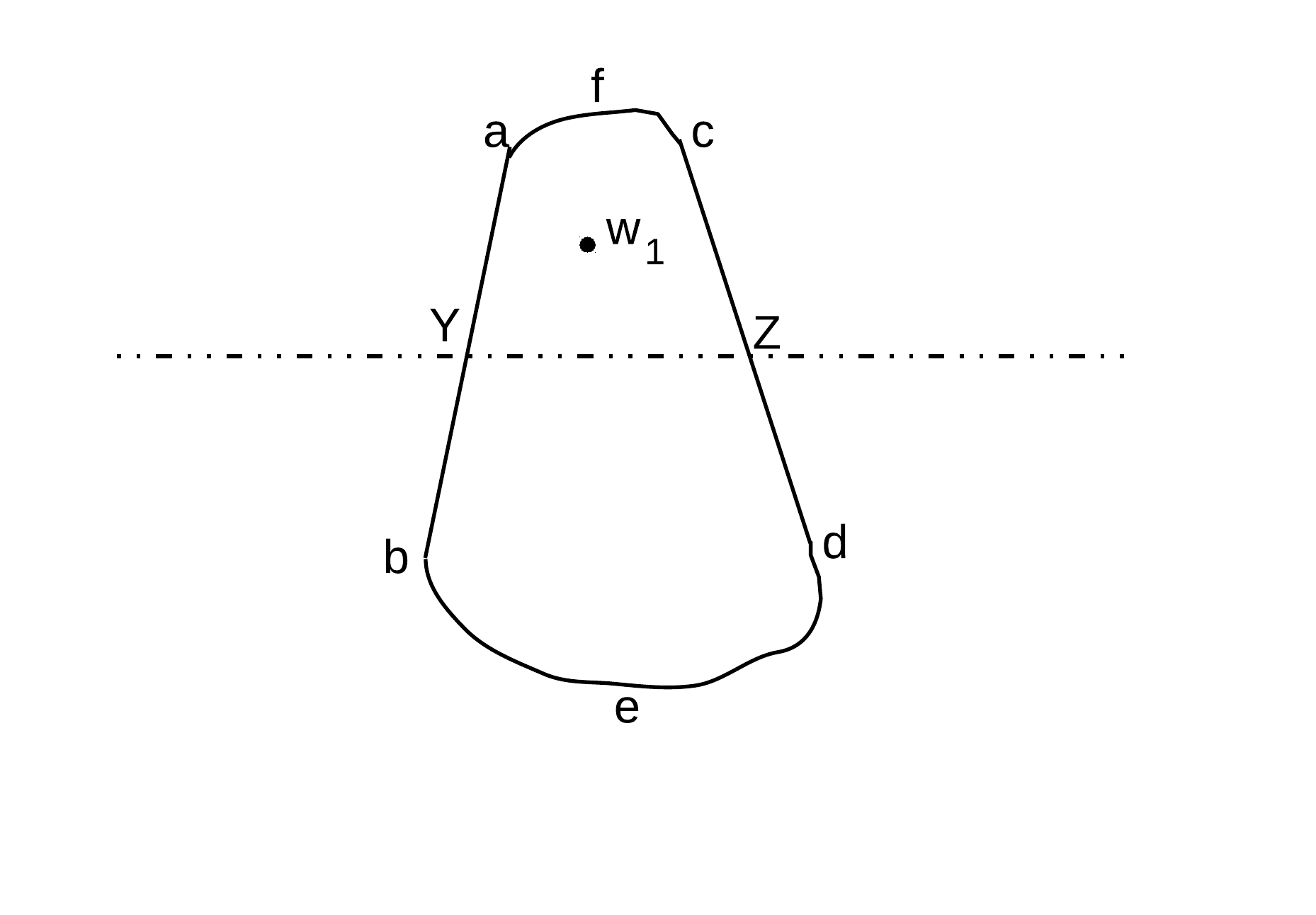}
        \caption{}
    \end{subfigure}

\caption{\((a)\) If~\(J_1\) lies below~\(\Gamma,\) then the path~\(\Gamma\) would cross the top side of~\(R.\) \((b)\) The top cell~\(T = afcdeba\) containing the dual vertex~\(w_1\) in its interior. The dual edge~\((w_1,w_2)\) with~\(w_2\) present in the interior of~\(R,\) must cross the top edge of~\(R\) in the segment~\(YZ.\)}
\label{w1_lies_below}
\end{figure}

Suppose~\(TD^+(R,V)\) also occurs and let~\(L = (J_1,\ldots,J_m)\) be a vacant plus connected top bottom crossing where~\(J_1\) is a top cell,~\(J_m\) is a bottom cell and every other~\(J_i\) is an interior cell. Let~\(w_i\) be the vertex of~\(G_d\) present in the cell~\(J_i\) so that~\(P := (w_1,w_2,\ldots,w_m)\) is a path in~\(G_d.\) We claim that the vertex~\(w_1\) necessarily above~\(\Gamma.\) This is because if~\(w_1\) were to be present below~\(\Gamma,\) then the cell~\(J_1 \in G\) containing~\(w_1\) in its interior also lies below~\(\Gamma.\) Let~\(L_{w_1}\) be the infinite line parallel to the left side of~\(R\) and containing the point~\(w_1.\) Since~\(w_1\) is contained in the interior of~\(J_1,\) the some edge of~\(J_1\) contains a point~\(u \in L_{w_1}\) lying above~\(w_1.\) Since~\(\Gamma\) lies above~\(J_1,\) there exists~\(v \in \Gamma\) lying above~\(u\) (see Figure~\ref{w1_lies_below}~\((a)\) for illustration). This would imply that some edge of~\(\Gamma\) crosses the top side of~\(R,\) a contradiction.

From the above paragraph we therefore get that the dual vertex~\(w_1\) is necessarily above~\(\Gamma\) and an analogous analysis implies that~\(w_m\) is below~\(\Gamma.\) Next, we argue that the ``first" edge~\((w_1,w_2) \in P\) is either present in the interior of~\(R\) or crosses the top edge of~\(R.\)
For illustration we consider a magnified top cell~\(T = afcdeba\) in Figure~\ref{w1_lies_below}\((b)\) containing edges~\(ab\) and~\(cd\) that intersect the top edge of~\(R\) (the dotted-dashed line). Because of the interior edge property~\((a2),\) the dual edge~\((w_1,w_2)\) must cross the top edge of~\(R\) in the segment~\(YZ.\)

Summarizing, the edge~\((w_1,w_2)\) is either present in the interior of~\(R\) or crosses the top edge of~\(R.\) Similarly the final edge~\((w_{m-1},w_m)\) is either contained in the interior of~\(R\) or crosses the bottom edge of~\(R.\) Every other vertex~\(w_i, 2 \leq i \leq m-1\) is present in the interior of~\(\Gamma.\) Therefore the path~\(P\) necessarily crosses~\(\Gamma\) in the sense that there are edges~\(f \in \Gamma\) and~\(e = (w_i,w_{i+1}) \in P\) such~\(f\) intersects~\(e.\) The end-vertices of the dual edge~\(e\) belong to cells~\(J_{i}\) and~\(J_{i+1}\) and the edge~\(f\) is common to~\(J_i\) and~\(J_{i+1}.\) At least one of the two cells~\(J_i\) or~\(J_{i+1}\) lies in the interior of~\(R\) and so the edge~\(f\) is necessarily contained within the rectangle~\(R.\) Thus~\(f\) is an edge of some occupied cell in the left right crossing~\(L_s\) and so one of the cells~\(J_i\) or~\(J_{i+1}\) must lie in the interior of~\(R\) and also be occupied, a contradiction.

\begin{figure}[tbp]
\centering
\includegraphics[width=3in, trim = 40 220 20 150, clip = true]{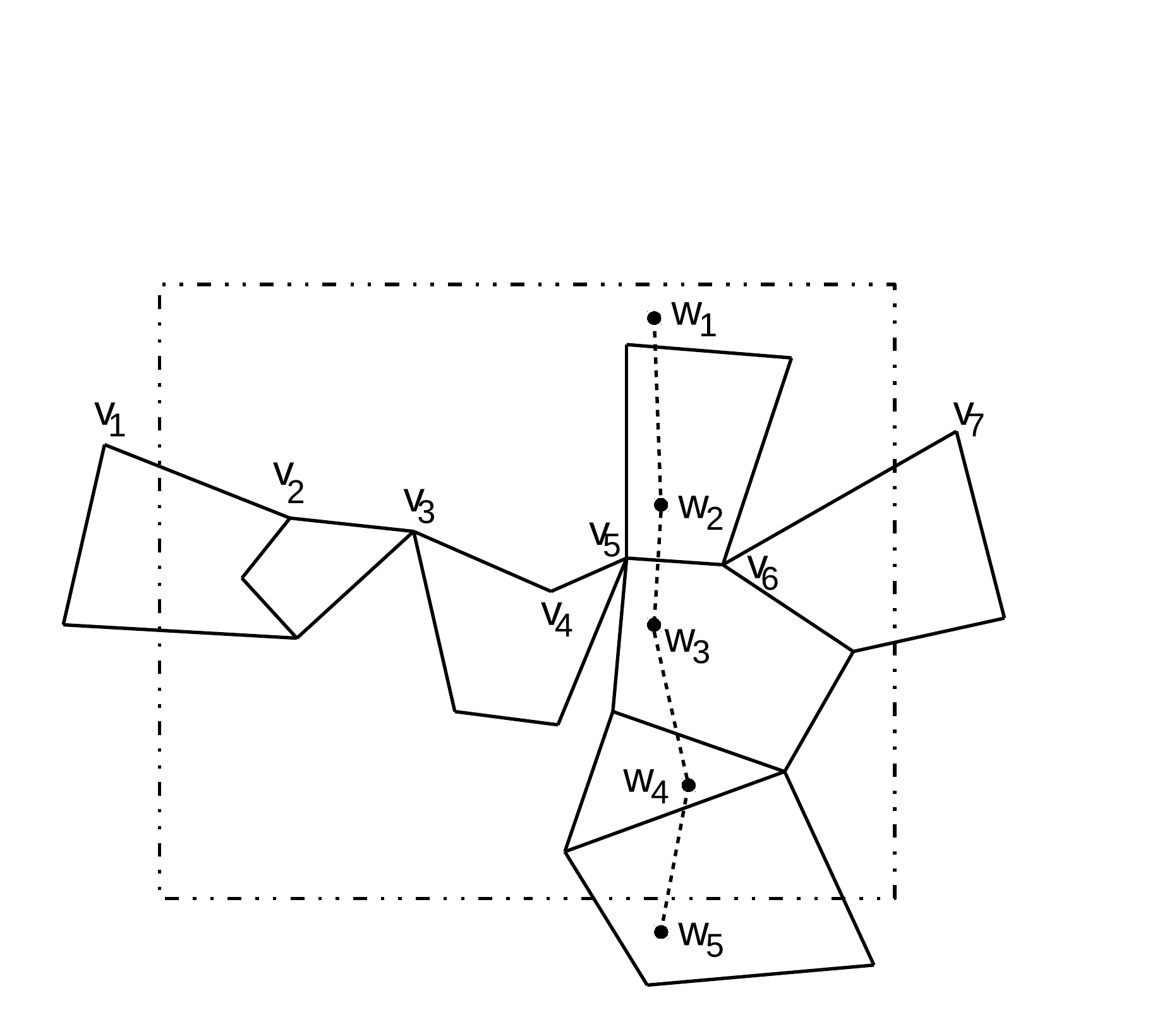}
\caption{The edge~\(f = (v_5,v_6)\) belonging to the path~\(\Gamma = (v_1,v_2,\ldots,v_7)\) and the dual edge~\(e = (w_2,w_3)\) belonging to the path~\(P = (w_1,\ldots,w_5)\) intersect.}
\label{lr_crossing_example}
\end{figure}

We illustrate the above argument in Figure~\ref{lr_crossing_example}, where the edge~\(f = (v_5,v_6) \) belonging to the path~\(\Gamma =(v_1,v_2,v_3,v_4,v_5,v_6)\) and the dual edge~\(e = (w_2,w_3)\) in the path~\(P = (w_1,w_2,w_3,w_4,w_5)\) intersect. ~\(\qed\)

\emph{Proof of~\((II)\)}: The collection~\({\cal I}_{tot}\) of all left cells~\(L_1,\ldots,L_l\) and the corner cells~\(C_{TL}\) and~\(C_{BL}\) in Figure~\ref{nice_cover} is a plus connected component. To each cell in~\({\cal I}_{tot}\) we now assign a \emph{label}~\(\omega.\) If some occupied cell~\(Q\) in the interior of~\(R\) is connected to a left cell by a star connected occupied path, we assign the label~\(\omega\) to~\(Q\) as well. The collection of all cells with the label~\(\omega\) is a star connected component which we denote as~\({\cal F}_{tot}.\) If~\(R\) is any cell star adjacent to some cell in~\({\cal F}_{tot}\) and not present in~\({\cal F}_{tot},\) we assign the label~\(\delta\) to~\(R.\)

By assumption, any vertex of a cell~\(Q \in {\cal F}_{tot}\) is contained in the interior of some dual cell. Therefore by Theorem~\ref{thm4}, there exists a plus connected cell cycle~\(\Delta_{vac} = (Z_1,\ldots,Z_t)\) surrounding~\({\cal F}_{tot}\) in such a way that the outermost boundary cycle~\(\partial_Z\) of~\(\Delta_{vac}\) contains all the cells of~\({\cal F}_{tot}\) in its interior. The cell cycle~\(\Delta_{vac}\) contains a plus connected cell sub-path~\(\Delta_Z = (Z_{u_1},\ldots,Z_{u_2})\) that lies to the right of~\({\cal L}_{left},\) with~\(Z_{u_1}\) intersecting~\({\cal L}_{top}\) and~\(Z_{u_2}\) intersecting~\({\cal L}_{bottom}.\)

\begin{figure}[tbp]
\centering
\includegraphics[width=3in, trim = 40 200 20 150, clip = true]{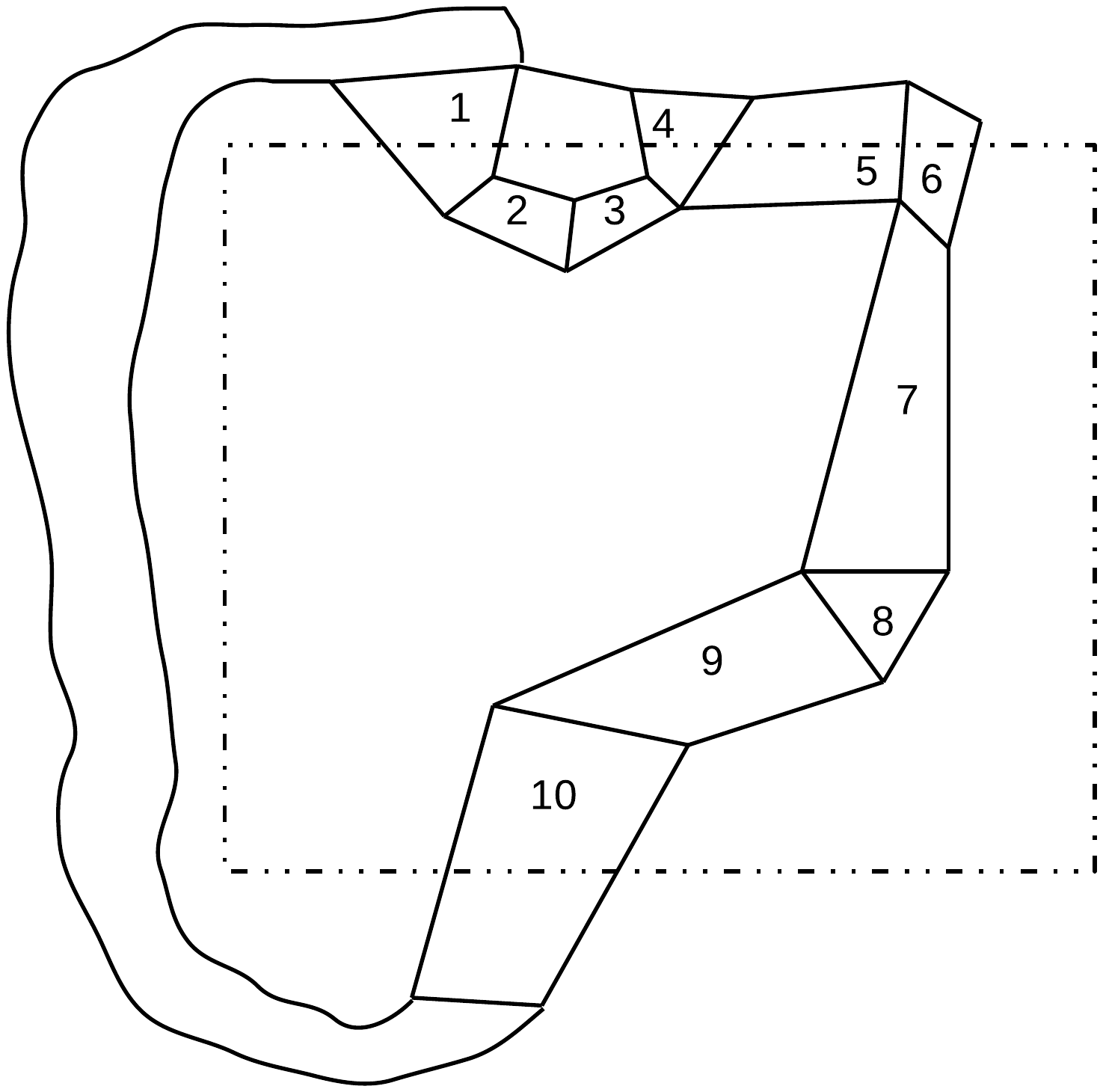}
\caption{The part of the vacant cycle~\(\Delta_{vac}\) that intersects the rectangle~\(R.\)}
\label{lr_contained_fig}
\end{figure}

In Figure~\ref{lr_contained_fig}, we illustrate the part of the cycle~\(\Delta_{vac}\) that intersects the rectangle~\(R\) with the cell labelled~\(i\) denoting~\(Z_i.\) As in Figure~\ref{lr_contained_fig}, the cycle~\(\Delta_{vac}\) may intersect the top side of~\(R\) multiple times but there exists a ``last" cell after which the cycle never intersects the top side of~\(R.\) Formally,~\(u_1\) be the largest index~\(i\) such that the cell~\(Z_{i}\) intersects the top side of~\(R\) and~\(Z_{u_2}\) is the ``first" cell after~\(Z_{u_1}\) that intersects the bottom side of~\(R.\) In Figure~\ref{lr_contained_fig},~\(u_1 = 6\) and~\(u_2 = 10.\)

By definition the cell~\(Z_{u_1}\) intersects the line~\({\cal L}_{top},\) the cell~\(Z_{u_2}\) intersects the line~\({\cal L}_{bottom}\) and that every other~\(Z_{j}, u_1 < j < u_2\) neither intersects~\({\cal L}_{top}\) nor intersects~\({\cal L}_{bottom}\) but lies between these two lines.
By the nicely padded assumption we have that~\(u_2 > u_1.\) No cell~\(Z_j, u_1 < j < u_2\) can be a right cell because then there would exist an occupied cell contained in the interior of~\(R\) which is star adjacent to~\(Z_j\) and is connected to some left cell~\(L_x\) by an occupied star connected cell path~\(P.\) The concatenation~\((L_x,P,Z_j)\) would then form an occupied star connected left right crossing of~\(R,\) a contradiction. Thus each cell~\(Z_j, u_1 < j < u_2\) is necessarily an interior cell.

By the nicely covered assumption, this necessarily means that~\(Z_{u_1}\) must be a top cell and not a corner cell. This is because, no corner cell is plus adjacent to an interior cell. Similarly~\(Z_{u_2}\) must be a bottom cell and not a corner cell and so~\((Z_{u_1},Z_{u_1+1},\ldots,Z_{u_2})\) forms a vacant plus connected top down crossing of~\(R.\)~\(\qed\)

\emph{Proof of~\((III)\)}: The proof is analogous to the proof of~\((II)\) with minor modifications. Here~\(\Delta_{vac} = (Z_1,\ldots,Z_t)\) is star connected and if the corner cell~\(C_{TR}\) or the bottom cell~\(C_{BR}\) in Figure~\ref{nice_cover} appear in~\(\Delta_{vac},\) we simply remove the corresponding entry from~\(\Delta_{vac}.\) The resulting sequence of vacant cells is still star connected and we proceed as before to get the desired vacant star connected top bottom crossing of~\(R.\)~\(\qed\)

\subsection*{\em Bond Percolation}
In this section, we consider bond percolation in the lattice~\(G\) and the mutual exclusivity of left right and top down crossings of in a rectangle.
We consider unoriented bond percolation and an analogous analysis holds for the oriented case as well.
As before we assume that the percolation lattice~\(G = \bigcup_{k=0}^{N}S_k\) and the dual lattice~\(G_d = \bigcup_{l=0}^{M}W_l\) satisfy properties~\((a1)-(a5)\) and the origin is the present in the interior of the cell~\(S_0.\)

Moreover, we also assume that for a fixed rectangle~\(R,\) we assume that the sides of~\(R\) nicely covered and nicely padded by cells of~\(G\) as described prior to Definition~\ref{lr_crossing_def}.






Assuming that~\(R\) is nicely covered as described above, we have the following definition of left right crossings.
\begin{Definition}\label{lr_crossing_def2}
A path~\(P = (v_1,\ldots,v_m) \subset G, m \geq 4\) is said to be a \emph{left right crossing} of a rectangle~\(R\) if:\\
\((d1)\) The edge~\((v_1,v_2)\) intersects the left side of~\(R\) and~\(v_1\) lies in the exterior of~\(R.\)\\
\((d2)\) The edge~\((v_{m-1},v_m)\) intersects the right side of~\(R\) and~\(v_m\) lies in the exterior of~\(R.\)\\
\((d3)\) Every other edge~\((v_i,v_{i+1}), 2 \leq i \leq m-2\) is contained in the interior of~\(R.\)
\end{Definition}
By definition any left right crossing must contain at least one edge in the interior of~\(R.\) Every edge in the closed interior of~\(R\) is now assigned one of the following two states: open or closed. Moreover, if edge~\(e \in G\) is open and if~\(f\) is the unique dual edge intersecting~\(G_d,\) then we assign~\(f\) to be open as well. If every interior edge in a left right crossing~\(P\) of~\(R\) is open, we say that \(P\) is an \emph{open left right crossing} of the rectangle \(R.\) An analogous definition holds for top down crossings. We denote \(LR\) to be the event that the rectangle~\(R\) contains an open left right crossing of~\(G.\)

For the dual crossing we have a slightly different definition. We say that~\(P_g := (g_1,\ldots,g_t)\) is a dual top bottom crossing of~\(R\) if the dual vertex~\(g_1\) lies in a top cell, the dual vertex~\(g_t\) lies in a bottom cell and each edge~\((g_i,g_{i+1}), 1 \leq i \leq t-1\) intersects an interior edge of~\(G;\) i.e., an edge of~\(G\) with both end-vertices present in the interior of~\(R.\)

We now see that every edge in a dual top bottom crossing has a state. By definition it suffices to see that the first edge~\((g_1,g_2)\) and the last edge~\((g_{t-1},g_t)\) both have states. First consider the edge~\((g_1,g_2)\) and let~\((v_1,v_2) \in G\) intersect~\((g_1,g_2).\) By the nicely covered assumption in Figure~\ref{nice_cover}, we see that the edge~\((v_1,v_2)\) belongs to one of the interior paths represented by the dotted lines and so necessarily lies in the interior of~\(R\) and consequently has a state. This implies that the dual edge~\((g_1,g_2)\) has the same state as~\((v_1,v_2).\)
Similarly, the last edge~\((g_{t-1},g_t)\) also has a state.  If every edge in~\(P_g\) is closed, we say that~\(P_g\) is a closed dual top bottom crossing and we let~\(TD_d\) be the event that~\(R\) contains a closed dual top bottom crossing consisting of edges of~\(G_d.\)


We have the following result.
\begin{Theorem}\label{thm_lr_bond} Suppose properties~\((a1)-(a5)\) hold. Further suppose that the rectangle~\(R\) is nicely covered and nicely padded by cells in~\(G\) as in Figure~\ref{nice_cover}. One of the events~\(LR\) or~\(TD_d\) always occurs, but not both.
\end{Theorem}



\begin{figure}[tbp]
\centering
\includegraphics[width=3in, trim = 5 160 10 140, clip = true]{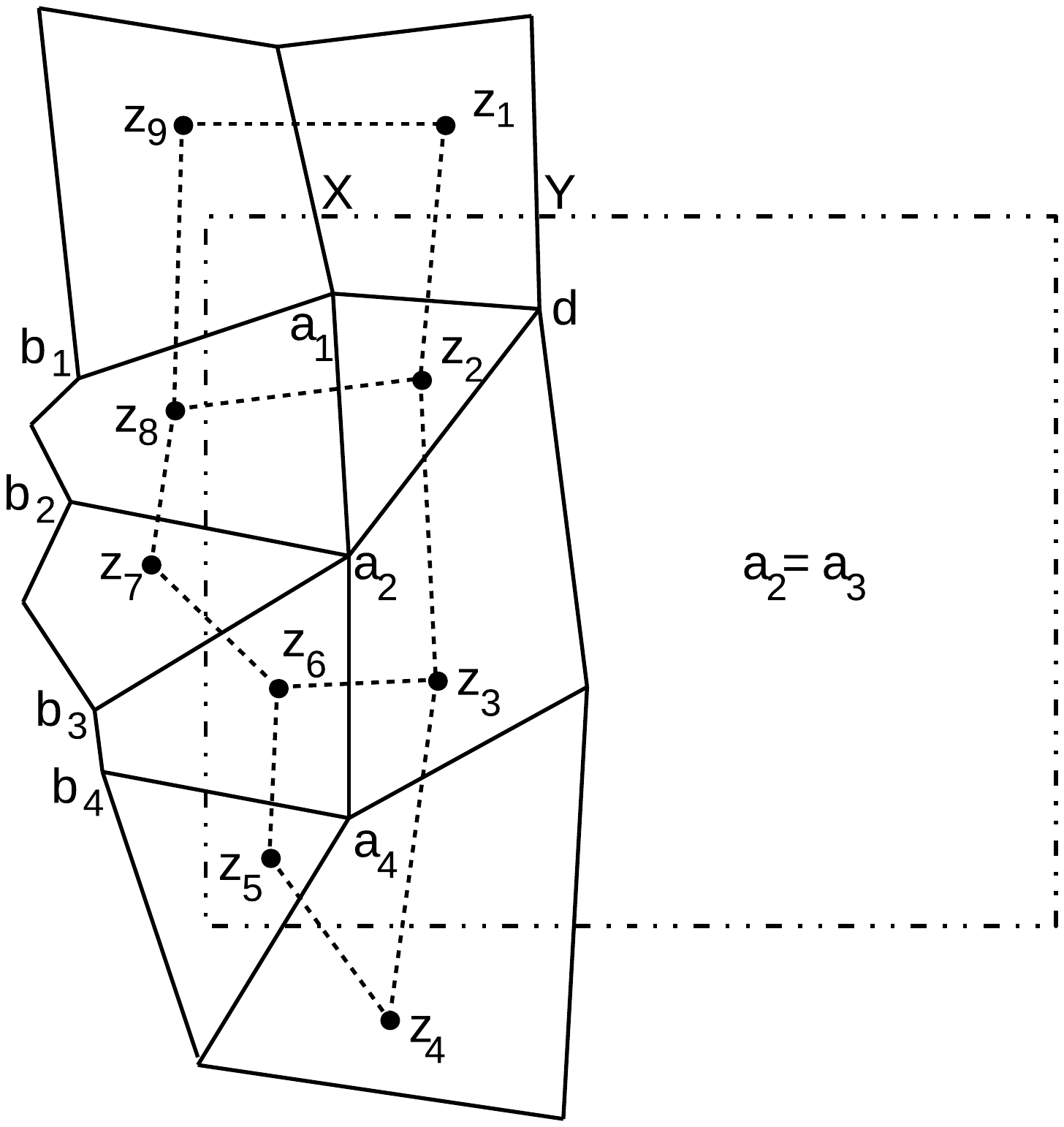}
\caption{The outermost boundary cycle~\(\partial_E = (z_1,\ldots,z_9,z_1)\) denoted by the dotted lines, intersects the top and bottom sides of the rectangle~\(R.\)}
\label{strong_covering}
\end{figure}


As before we need to prove three statements:\\
\((I)\) Both~\(LR\) and~\(TD_d\) cannot occur simultaneously.\\
\((II)\) If~\(LR\) does not occur, then~\(TD_d\) occurs.\\
\((III)\) If~\(TD_d\) does not occur, then~\(LR\) occurs.\\
The proof of~\((I)\) is analogous to the proof of~\((I)\) in Theorem~\ref{thm_lr}. If~\(\Gamma\) is any open left right crossing and~\(\Delta\) is any top bottom dual crossing, we obtain that one vertex of the dual crossing lies above~\(\Gamma\) and one vertex lies below~\(\Gamma\) and so these two paths must necessarily meet.  The dual edge intersecting any edge~\(e \in G\) has a state and in fact the same state as~\(e\) and this leads to a contradiction.~\(\qed\)


The proof of~\((III)\) is analogous to~\((II)\) and we prove~\((II)\) below.\\
\emph{Proof of~\((II)\)}: Let~\(\{e_i\}_{1 \leq i \leq t}\) be the set of edges of~\(G\) intersecting the left side of~\(R\) arranged in the decreasing order of the~\(y-\)coordinate of the intersection point and for~\(1 \leq i \leq t,\) let~\(a_i\) and~\(b_i\) be the end-vertices of~\(e_i\) present in the interior and exterior of~\(R,\) respectively.  For example, in Figure~\ref{nice_cover}, the edge~\(e_1=a_1b_1, e_2 = a_2b_2\) and so on.

Let~\({\cal I}_{tot}\) be the set of all open edges lying in the interior of~\(R\) and connected to some vertex in~\(\{a_i\}_{1 \leq i \leq t}\) by an open path and for~\(1 \leq i \leq t-1,\) let~\(A_i\) be the path between~\(a_i\) and~\(a_{i+1}\) lying in the interior of the rectangle~\(R.\) The union~\({\cal E}_{tot} = {\cal I}_{tot} \cup \{A_i\}_{1 \leq i \leq t-1}\) is then a connected component and  each vertex~\(v \in {\cal E}_{tot}\) is present in the interior of some dual cell~\(W(v) \subset G_d.\) The union of the dual cells~\(\{W(v)\}_{v \in {\cal E}_{tot}}\) forms a plus connected dual component whose outermost boundary~\(\partial_W\) is a single cycle in~\(G_d\) containing all edges of~\({\cal E}_{tot}\) in its interior.


We now see that~\(\partial_W\) contains at least one dual vertex present in the interior of a top cell. From Figure~\ref{nice_cover}, the dual edge joining~\(x\) and~\(y\) intersects the edge~\((a_1,d)\) and so belongs to the dual cell~\(W(a_1)\) containing~\(a_1\) in its interior. The dual vertex~\(y\) lying in the interior of the top cell~\(T_1\) therefore belongs to the dual cell~\(W(a_1),\) by property~\((a5).\) The dual edge~\((x,y)\) is not present in any other dual cell~\(W(v), v \in {\cal E}_{tot}\) and so belongs to the final cycle~\(\partial_W\) as well.

By an analogous argument, all the dual vertices present in the interior of the left cells~\(L_1,\ldots, L_l\) and the corner cells~\(C_{BL}\) and~\(C_{TL}\) form a sub-path~\(P_W\) of~\(\partial_W\) and the dual edge~\((x_0,y_0) \in \partial_W\) as well, where~\(x_0\) and~\(y_0\) are the dual vertices is  present in the interior of the bottom left corner cell~\(C_{BL}\) and the first bottom cell~\(B_1,\) respectively (See Figure~\ref{nice_cover}). The subpath~\(\Delta_W := \partial_W \setminus P_W\) has end-vertices~\(y\) and~\(y_0.\) For illustration, in Figure~\ref{strong_covering}, the outermost boundary cycle~\(\partial_E\) formed by the merging of the dual cells~\(\{W(a_i)\}_{1 \leq i \leq t}\) is shown by the dotted line~\((z_1,z_2,\ldots,z_9,z_1).\) Here~\(z_1 = y,z_9= x,z_5 = x_0\) and~\(z_4 = y_0.\) The sub-path~\(P_W = (z_5,z_6,z_7,z_8,z_9)\) and~\(\Delta_W =(z_1,z_2,z_3,z_4).\)

The path~\(\Delta_W\) might contain many dual vertices present in the interior of some top cell and so we pick the ``last" such vertex and call it~\(y_1.\)
Similarly we pick the first dual vertex present in the interior of a bottom cell. Formally, we extract a sub-path~\(\Gamma_W := (y_1,\ldots,y_p) \subset \Delta_W\) such that the vertex~\(y_1\) lies in the interior of a top cell, the vertex~\(y_p\) lies in the interior of a bottom cell and every other~\(y_i\) lies in the interior of a right, corner or interior cell. To prove that each edge of~\(\Gamma_W\) has a state it suffices to see that there does not exist~\(m\) such that~\(y_m\) belongs to a corner or a right cell and~\(y_{m+1}\) belongs to a right cell.

\begin{figure}[tbp]
\centering
\includegraphics[width=3in, trim = 40 220 20 130, clip = true]{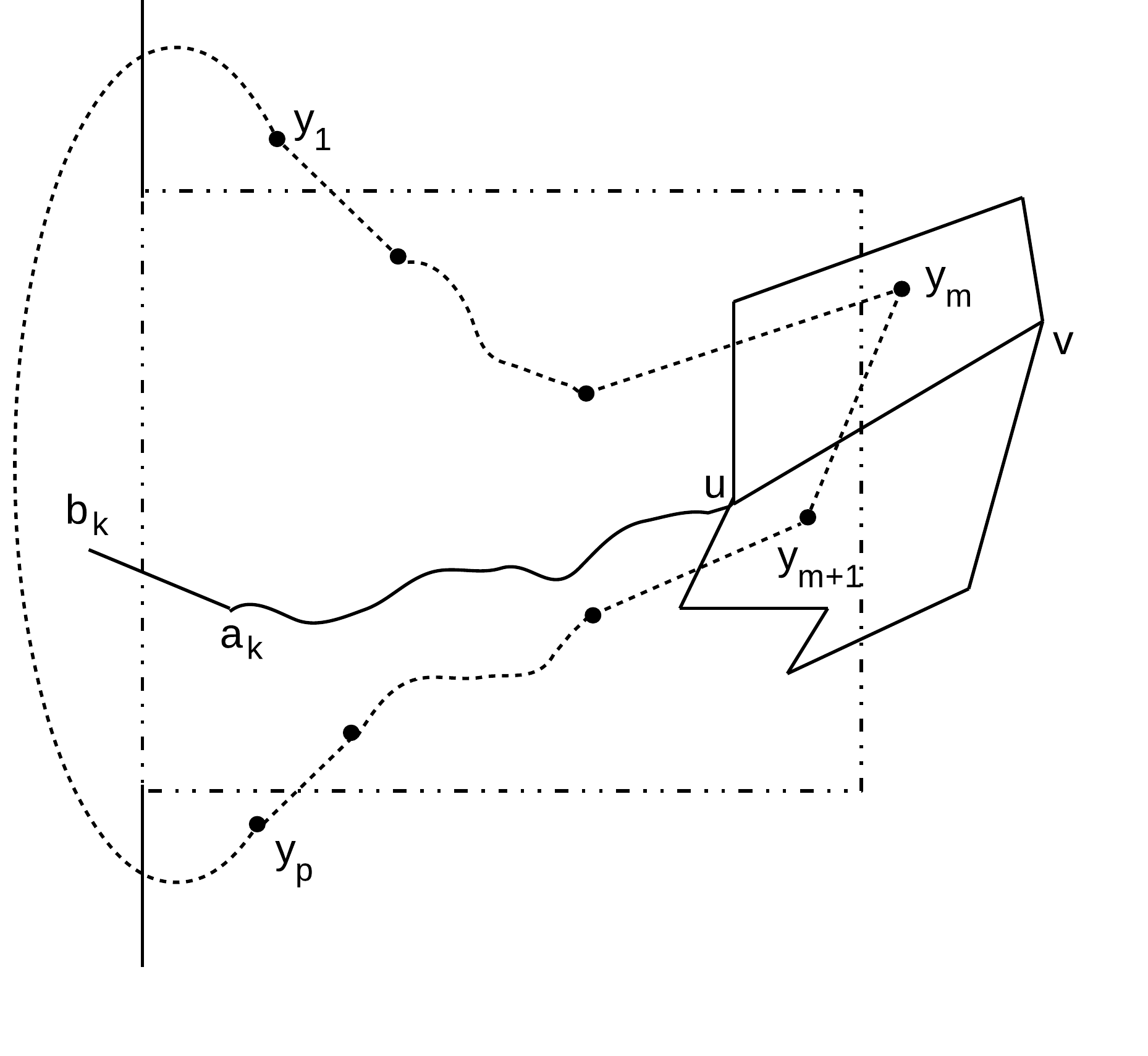}
\caption{The edge~\((y_m,y_{m+1}) \in \Gamma_W\) cuts an edge intersecting the right side of~\(R.\) }
\label{lr_crossing_out}
\end{figure}

As in Figure~\ref{lr_crossing_out}, we assume that~\(y_m\) and~\(y_{m+1}\) both belong to right cells and an analogous argument holds for the corner cells since no corner cell is plus adjacent with an interior cell. From Figure~\ref{lr_crossing_out} we see that the edge~\((y_{m},y_{m+1}) \in \Gamma_W\) intersects some edge~\((u,v)\) that cuts the right side of~\(R\) as in Figure~\ref{lr_crossing_out}. The vertex~\(u \in G\) is therefore contained in the interior of the dual cell containing~\((y_m,y_{m+1})\) (property~\((a5)\)) and so by definition of the component~\({\cal E}_{tot},\) there exists an open path~\(P\) from~\(u\) to some vertex~\(a_k\) of the edge~\(e_k\) that intersects the left side of~\(R.\) The concatenation~\((e_k,P,(u,v))\) would then form an open left right crossing of~\(R,\) a contradiction.

Finally by the nicely padded assumption there must exist at least  two edges in~\(\Gamma_W\) and so the subpath~\((y_1,\ldots,y_p) \subset \partial_W\)  is a dual top bottom crossing of~\(R.\) It remains to see that each such edge is closed. Suppose not and the edge~\((y_i,y_{i+1}) \in \Gamma_W\) is open. By the interior edge property~\((a2),\) there exists exactly one edge~\((v_1,v_2) \in G\) that intersects~\((y_i,y_{i+1}).\) By construction, one of the vertices, say~\(v_1\) belongs to the component~\({\cal E}_{tot}\) and so there is an open path from~\(v_1\) to some end-vertex~\(a_k\) of the edge~\(e_k\) intersecting the left side of~\(R.\)

Since~\((y_i,y_{i+1})\) is open, the edge~\((v_1,v_2)\) is open as well and so~\(v_2\) also belongs to~\({\cal E}_{tot}.\) But if~\(R(v_1)\) and~\(R(v_2)\) denote the dual cells containing~\(v_1\) and~\(v_2,\) respectively, then from property~\((a4)\) the cells~\(R(v_1)\) and~\(R(v_2)\) are plus adjacent and share the edge~\((y_i,y_{i+1}).\) This implies that~\((y_i,y_{i+1})\) is present in the interior of the cycle formed by merging~\(R(v_1)\) and~\(R(v_2)\) and consequently~\((y_i,y_{i+1})\) must be present in the interior of the outermost boundary cycle~\(\partial_W\) as well, a contradiction.~\(\qed\)

\begin{acknowledgments}
I thank Professors Rahul Roy, Thomas Mountford, Federico Camia, C. R. Subramanian and the referee for crucial comments that led to an improvement of the paper. I also thank IMSc for my fellowships.
\end{acknowledgments}

\bibliographystyle{plain}


\end{document}